\documentclass[]{article}

\title{Existence, multiplicity and stability of endemic states for an 
 age-structured S--I epidemic model}

\author{D.~Breda\thanks{Department of Mathematics and Informatics, 
   University of Udine, via delle Scienze 208, 33100 Udine, Italy.}\; and\ 
 D.~Visetti\thanks{Department of Mathematics, University of Trento, 
   via Sommarive 14, 38123 Povo (TN), Italy.}}

\RequirePackage[italian,english]{babel} 
\RequirePackage[dvips]{graphicx}
\RequirePackage{amsmath}
\RequirePackage{amssymb}
\RequirePackage{amsthm}

\numberwithin{equation}{section}

\begin{document}
\selectlanguage{english}

\newcommand{\NN}{{\mathbb N}}
\newcommand{\RR}{{\mathbb R}}
\newcommand{\ZZ}{{\mathbb Z}}
\newcommand{\QQ}{{\mathbb Q}}
\newcommand{\CC}{{\mathbb C}}

\newtheorem{theorem}{Theorem}[section]
\newtheorem{proposition}[theorem]{Proposition}
\newtheorem{lemma}[theorem]{Lemma}
\newtheorem{corollary}[theorem]{Corollary}
\newtheorem{remark}[theorem]{Remark}
\newtheorem{definition}[theorem]{Definition}


\maketitle

\begin{abstract}
 We study an S--I type epidemic model in an age-structured 
 population, with mortality due to the disease. A threshold 
 quantity is found that controls the stability of the 
 disease-free equilibrium and guarantees the existence of an 
 endemic equilibrium. We obtain conditions on the age-dependence of the
 susceptibility to infection that imply the 
 uniqueness of the endemic equilibrium. An example with two 
 endemic equilibria is shown.  Finally, we analyse numerically 
 how the stability of the endemic equilibrium is affected by 
 extra-mortality and by the possible periodicities induced by 
 the demographic age-structure.
\end{abstract}


\section{Introduction}

In this paper we study the existence and uniqueness or multiplicity of 
positive steady states for an S--I type epidemic through 
an age-structured population.  The population growth is subject to 
intra-specific mechanisms that, in the absence of infection, lead to a 
steady state which can be destabilised by Hopf bifurcation, for some 
values of the intrinsic reproduction number $R_0^d$. The infection is 
regulated by an age-dependent contact mechanism and interferes with 
the demographic process introducing additional mortality due to the disease.

A model with disease-induced mortality was for the first time introduced in  
\cite{AM} and the formulation in an age-structured context is due to May in 
\cite{M}.   In \cite{Andreasen} Andreasen 
analysed two different cases:  an age-structured epidemiological model and 
a model where infection has a constant duration.   For the former he assumed 
that only the epidemiological interactions are density-dependent;  
he obtained conditions for the existence of one endemic state and 
studied its stability in the hypothesis that only fecundity is age dependent.

Here we give conditions for existence and uniqueness of an endemic 
equilibrium.   A situation where two endemic equilibria appear is shown.  
Multiple solutions are specific of this model, in the sense that they do not
occur in models without age structure or in the same model without the
extra-mortality. Models with multiple steady states had already been 
proposed in the literature, the current model exhibits a different 
mechanism through which they arise.  In particular, in \cite{huang} the 
role of variable population size in a multigroup epidemic model is 
emphasised, showing that multiple endemic equilibria are possible.   
Concerning the analysis of epidemic models with multiple endemic 
equilibria in an SIS epidemic model without age structure, 
see also \cite{limazhou} and the bibliography therein.

In addition to the study of existence of equilibria, we also investigate 
their stability showing that the effect of disease-induced mortality 
may produce stability changes related to Hopf's 
bifurcation. Using a numerical method exposed in \cite{breda}, we 
explore the behaviour of our model, within some critical parameter regions.

The paper is structured as follows.  Next section is devoted to present our 
model, while in Section \ref{sec:wp} the well-posedness of the 
system of partial differential equations is pointed out.  In 
Section \ref{sec:endeq} conditions for existence of 
endemic equilibria are found.  Successively, in Section \ref{secexun} 
uniqueness versus multiplicity of endemic 
states is discussed.

In Section \ref{sec:ce} the characteristic equation is 
computed, so that in Section \ref{sec:stabdfe} the stability analysis 
of the disease-free equilibrium is carried out and in Section 
\ref{secsette} a stability change is described.

Finally, in Section \ref{sec:stab_an} stability is numerically analysed 
in two different situations:  where multiple endemic states occur and where 
periodic solutions arise through Hopf's bifurcation.  Our analysis reveals 
the importance of the parameter of disease-induced mortality.   In the 
last section conclusions are drawn.

\section{The model}

We consider a population that, in the absence of infection is described 
by  $N(a,t)$, the age density at time $t \ge 0$, where $a\in [0,a^\dagger]$, 
$a^\dagger < +\infty$ being the maximum age an individual of the population 
may reach. The growth of the population is regulated by the following 
model of  Gurtin-MacCamy type (see \cite{gurtin})
\begin{equation}
\label{Gurtin_MC}
\left\{ \begin{aligned}
   &\frac{\partial N}{\partial t}+
       \frac{\partial N}{\partial a}+\mu(a)N(a,t)=0\, , \\
   &N(0,t)= R_0^d \Phi(Q(t))\int_0^{a^\dagger}\beta(a)N(a,t)\,da\, .
 \end{aligned}\right.
\end{equation}

In (\ref{Gurtin_MC}), $\beta(a)$ and $\mu(a)$ are the intrinsic vital 
rates, with $\beta$ normalised to satisfy
\begin{equation*}
 \label{intbeta}
 \int_0^{a^\dagger}\beta(\sigma)\pi(\sigma)\, d\sigma = 1\, ,
\end{equation*}
where
\begin{equation*}
 \label{pi}
 \pi(a)=e^{-\int_0^a\mu(\sigma)\,d\sigma}
\end{equation*}
is the survival probability, i.e. the probability at birth 
of surviving to age $a$. Moreover, the parameter $R_0^d$ is the 
demographic basic reproduction ratio and we suppose that
$$
R_0^d>1\, .
$$
The function $\Phi$ is a Lipschitz-continuous non-increasing 
function describing density-dependence of births.  We assume that
$$
\Phi(0)=1\, ,\qquad \lim_{x\to +\infty}\Phi(x)=0
$$
and that there exists $x_0\in\RR$ such that $\Phi(x_0)<1/R_0^d$ 
and $\Phi(\cdot)$ is decreasing on $[0,x_0]$.  Finally, $Q$ is the 
size
\begin{equation}\label{qu}
Q(t) = \int_0^{a^\dagger}r(\sigma)N(\sigma,t)\, d\sigma ,
\end{equation}
where $r(a)$ is a weight kernel. We point out that, by the 
assumptions made, equation
\begin{equation}
 \label{Qd}
 R_0^d\Phi(Q_d)=1
\end{equation}
has a unique solution $Q_d^*$, providing a non trivial equilibrium
\begin{equation}
 \label{N*}
 N^*(a)=\frac{Q_d^*\pi(a)}{\int_0^{a^\dagger}r(\sigma)\,\pi(\sigma)\, d\sigma}
\end{equation}
of the population.

We recall (see \cite{Iannelli}) that, provided that $\Phi$ is differentiable at $Q^*_d$, 
the stability of this steady state is related to the characteristic equation
\begin{equation}\label{carN}
1= \int_0^{a^\dagger } e^{-\lambda a}\beta(a)\pi(a) da + R_0^d 
 \phi'(Q_d^*) N^*(0)  \int_0^{a^\dagger } e^{-\lambda a}r(a)\pi(a) da\, ,
\end{equation}
the analysis of which also provides information about possible 
bifurcation points.

For the problem above we make the following standard hypotheses:
$$
\begin{aligned}
\beta(\cdot)\in L^\infty(0,a^\dagger),\hspace{1cm}&
\beta(a)\geq 0\mbox{ in } [0,a^\dagger],\\
\mu(\cdot)\in L^1_{\mbox{loc}}(0,a^\dagger),\hspace{1cm}&
\mu(a)\geq 0\mbox{ in } [0,a^\dagger],\\
r(\cdot)\in L^\infty(0,a^\dagger), \hspace{1cm}&r(a)\geq 0\mbox{ in }
[0,a^\dagger]\, .
\end{aligned}
$$  Since no individual may live past age 
$a^\dagger$, in order to have
\[
\pi(a^\dagger) =0,
\]
we need to assume
$$
\int_0^{a^\dagger}\mu(\sigma)\, d\sigma=+\infty\, .
$$

The population develops an S--I type epidemics so that, 
denoting respectively  by $S(a,t)$ and $I(a,t)$ the age-specific densities of susceptible and 
infective individuals at time $t$,  the following system of equations 
describes the transmission dynamics of the disease:
\begin{equation}
 \label{eq}
 \left\{ \begin{aligned}
     &\frac{\partial S}{\partial t}+
       \frac{\partial S}{\partial a}=-\left(\lambda(a,t)+\mu(a)\right)
       S(a,t)\, , \\
     &\frac{\partial I}{\partial t}+
       \frac{\partial I}{\partial a}=\lambda(a,t)\,S(a,t)-\left(\mu(a)
       +\alpha\right)I(a,t)\, , \\
     &S(0,t)=R_0^d\Phi(Q(t))\int_0^{a^\dagger}\beta(a)\left(
       S(a,t)+I(a,t)\right)da\, , \\
     &I(0,t)=0\, , \\
     &S(a,0)=S_0(a)\geq 0\, , \\
     &I(a,0)=I_0(a)\geq 0\, .
 \end{aligned}\right.
\end{equation}
Here the terms $\lambda(a,t)$ and $\alpha$ describe the mechanism of 
infection.   In particular, $\alpha$ is the extra-mortality due to the 
infection and $\lambda(a,t)$ is the force of infection.  Concerning 
$Q(t)$, in (\ref{eq}) we set
\begin{equation}
 \label{Q}
 Q(t) = \int_0^{a^\dagger}r(\sigma)(S(\sigma,t)+I(\sigma,t))\, d\sigma ,
\end{equation}
where we implicitly assume (compare with (\ref{qu})) that both 
susceptible and infective individuals are equally active in the 
population. In fact, in (\ref{eq}) we are also assuming that newborns are all susceptible and 
equally produced  by both kind of individuals. We note that the 
variable
\[
N(a,t)=S(a,t)+I(a,t)
\] 
satisfies the equation
$$
\frac{\partial N}{\partial t}+
 \frac{\partial N}{\partial a}+\mu(a)N(a,t)+\alpha I(a,t)=0\, ,
$$
showing how the disease-induced mortality affects the population 
growth.

For the parameters regulating the infection mechanism we assume that 
$\alpha\geq 0$ and consider the general separable inter-cohort form 
of the force of infection (see \cite{bci}):
\begin{equation*}
 \lambda(a,t) = K(a)\int_0^{a^\dagger}q(\sigma)\, I(\sigma,t)\, 
 d\sigma\, ,
\end{equation*}
where the age-specific infectiousness $q(\cdot)$ and the age-specific 
contagion rate $K(\cdot)$ satisfy the following conditions
\[
q(\cdot),\, K(\cdot)\in L^\infty(0,a^\dagger),
\]
\[
q(a),\,  K(a)\geq 0\mbox{ in }[0,a^\dagger]\, .
\]

\section{Well-posedness}
\label{sec:wp}

Existence and uniqueness of a solution to problem (\ref{eq}) can be 
proved following the standard techniques presented for example 
in \cite{Iannelli} (see also \cite{tesi}). For the sake of completeness 
here we give a sketch of the proof without going into technical details.
First we denote
\begin{equation}
 \label{B(t)}
 B(t)=S(0,t)\, ,
\end{equation}
so that, integrating $S$ along the characteristics, we obtain
\begin{equation}
 \label{Scar}
 S(a,t)=\left\{\begin{array}{ll}
     S_0(a-t)e^{-\int_0^t[\lambda(a-t+\tau,\tau)+\mu(a-t+\tau)]\, d\tau} 
     & \mbox{for }a\geq t\, , \\
     B(t-a)\pi(a)e^{-\int_0^a\lambda(\sigma,t-a+\sigma)\, d\sigma}
     & \mbox{for }t\geq a\, .
   \end{array}\right.
\end{equation}
Moreover, using (\ref{Scar}) and integrating $I$ along the characteristics, 
we have
\begin{equation}
 \label{Icar}
 \begin{array}{l}
  \hspace{-2mm}I(a,t)=\\
   \hspace{-2mm}\left\{\begin{array}{ll}
     \hspace{-2mm}e^{-\alpha t-\int_0^t\mu(a-t+\tau)\, d\tau}\left(I_0(a-t)
       \phantom{\int_0^t}\right. & \\
     \hspace{-2mm}\quad\left.+S_0(a-t)
       \int_0^t\lambda(a-t+\tau,\tau)e^{\alpha\tau-\int_0^\tau\lambda
       (a-t+\xi,\xi)\, d\xi}d\tau\right)
     & \hspace{-3mm}\mbox{for }a\geq t\, , \\
     \hspace{-2mm}\pi(a)e^{-\alpha a}B(t-a)
       \int_0^a\lambda(\sigma,t-a+\sigma)e^{\alpha\sigma-\int_0^\sigma
         \lambda(\rho,t-a+\rho)\, d\rho}d\sigma
     & \hspace{-3mm}\mbox{for }t\geq a\, .
   \end{array}\right.
 \end{array}
\end{equation}
Then, if we denote
\begin{equation}
 \label{W(t)}
 W(t)=\int_0^{a^\dagger}q(a)I(a,t)\, da\, ,
\end{equation}
we obtain that $B(t)$ satisfies the following integral equation (see 
(\ref{eq})):
\begin{equation}
\label{B(t)int}
\begin{aligned}
B(t) = &R_0^d\Phi(Q(t)) \left\{ \int_0^t \beta(a)\pi(a) \Bigg( 
       e^{-\alpha a} \right. \\
     &\left. +\alpha\int_0^a e^{-\alpha(a-\sigma)-
       \int_0^\sigma K(\rho)W(t-a+\rho)\, d\rho}\, d\sigma 
       \right) B(t-a)\,da \\
     &+ \int_t^{+\infty} \beta(a) \Bigg[ S_0(a-t) \Bigg( e^{-\alpha t
       -\int_0^t \mu(a-t+\sigma)\, d\sigma} \\
     &+\alpha e^{-\int_0^t \mu(a-t+\sigma)\, d\sigma} \int_0^t 
       e^{-\alpha(t-\sigma)-\int_0^\sigma K(a-t+\tau)\, W(\tau)\, d\tau}
       d\sigma \Bigg) \\
     &+ I_0(a-t) e^{-\alpha t-\int_0^t \mu(a-t+\sigma)\, 
       d\sigma} \Bigg] da \Bigg\}\, ,
\end{aligned}
\end{equation}
where all the functions are extended by zero outside the interval $[0,
a^\dagger]$ and $Q$ is defined in (\ref{Q}).   To solve our problem, we 
first consider $Q$ and $W$ as two given continuous, nonnegative 
functions. It is well known that equation (\ref{B(t)int}), being a 
linear integral equation of Volterra type with a nonnegative kernel, 
admits a unique continuous and nonnegative 
solution.  We denote it by $B(t,Q,W)$ to show the dependence on $Q$ 
and $W$.

Now, for a fixed $T>0$, we consider the space $E=C\left([0,T];
(L^1(0,a^\dagger))^2\right)$ and the closed set
$$
\begin{aligned}
\mathcal{K} = &\{ k\in E\;|\; k(t)=(s(a,t),i(a,t)),\, 
   s(a,t)\geq 0,\, i(a,t)\geq 0,\\
 &\|s(\cdot,t)\|_{L^1(0,a^\dagger)}\leq M,\, \|i(\cdot,
   t)\|_{L^1(0,a^\dagger)}\leq M\}\, ,
\end{aligned}
$$
with
$$
M = (1+\alpha a^\dagger)e^{R_0^d\|\beta\|_{L^\infty(0,
 a^\dagger)}(1+\alpha a^\dagger)T}\| S_0+I_0\|_{L^1(0,
 a^\dagger)}\, .
$$
Then, given $k\in E$, $k(t)=(s(a,t),i(a,t))$, we set
\begin{eqnarray*}
Q(t) & = & \int_0^{a^\dagger} r(a)\,(s(a,t)+i(a,t))\, da\, , \\
W(t) & = & \int_0^{a^\dagger} q(a)\, i(a,t)\, da\, ,
\end{eqnarray*}
and define the map
$$
\mathcal{T}:\mathcal{K}\subset E\longrightarrow E\, ,
$$
where $\mathcal{T}(k)=(\tilde s(a,t),\tilde \imath(a,t))$ and
$$
\begin{aligned}
\tilde s(a,t) &= \left\{ \hspace{-2mm}\begin{array}{ll}
   S_0(a-t)e^{-\int_0^t[K(a-t+\sigma)W(\sigma)+\mu(a-t+\sigma)]\, d\sigma} 
   & \mbox{for }a\geq t\, , \\
   B(t-a,Q,W)\pi(a)e^{-\int_0^aK(\sigma)W(t-a+\sigma)\, d\sigma}
   & \mbox{for }t\geq a\, ,
 \end{array}\right. \\
\tilde\imath(a,t) &= \left\{ \hspace{-2mm}\begin{array}{ll}
   e^{-\alpha t-\int_0^t\mu(a-t+\sigma)\, d\sigma}\left(I_0(a-t)
     \phantom{\int_0^t}+S_0(a-t)\cdot\right. & \\
   \ \left.
     \cdot\int_0^tK(a-t+\sigma)W(\sigma)e^{\alpha\sigma-\int_0^\sigma 
       K(a-t+\tau)W(\tau)\, d\tau}d\sigma\right)
   & \mbox{for }a\geq t\, , \\
   \pi(a)e^{-\alpha a}B(t-a,Q,W)\cdot & \\
   \ \cdot\int_0^a K(\sigma)W(t-a+\sigma)
     e^{\alpha\sigma-\int_0^\sigma K(\tau)W(t-a+\tau)\, d\tau}d\sigma
   & \mbox{for }t\geq a\, .
 \end{array}\right.
\end{aligned}
$$

Standard estimates allow to prove that $\mathcal{T}$ sends 
$\mathcal{K}$ into itself and that a suitable power of $\mathcal{T}$ 
is a contraction, so that there exists a unique fixed point.  
Then, we can state the following result:

\begin{theorem}
 Let $(S_0,I_0)\in (L^1(0,a^\dagger))^2$, then there is one and 
 only one $k\in\mathcal{K}$, $k(t)=(S(a,t),I(a,t))$ such that
 $$
 \begin{aligned}
 S(a,t) &= \left\{ \hspace{-2mm}\begin{array}{ll}
     S_0(a-t)e^{-\int_0^t[K(a-t+\sigma)W(\sigma)+\mu(a-t+\sigma)]\, d\sigma} 
     & \mbox{for }a\geq t\, , \\
     B(t-a,Q,W)\pi(a)e^{-\int_0^aK(\sigma)W(t-a+\sigma)\, d\sigma}
     & \mbox{for }t\geq a\, ,
   \end{array}\right. \\
 I(a,t) &= \left\{ \hspace{-2mm}\begin{array}{ll}
     e^{-\alpha t-\int_0^t\mu(a-t+\sigma)\, d\sigma}\left(I_0(a-t)
       \phantom{\int_0^t}+S_0(a-t)\cdot\right. & \\
     \ \left.
       \cdot\int_0^tK(a-t+\sigma)W(\sigma)e^{\alpha\sigma-\int_0^\sigma 
       K(a-t+\tau)W(\tau)\, d\tau}d\sigma\right)
     & \mbox{for }a\geq t\, , \\
     \pi(a)e^{-\alpha a}B(t-a,Q,W)\cdot & \\
     \ \cdot\int_0^a K(\sigma)W(t-a+\sigma)
       e^{\alpha\sigma-\int_0^\sigma K(\tau)W(t-a+\tau)\, d\tau}d\sigma
     & \mbox{for }t\geq a\, ,
   \end{array}\right. \\
 Q(t) &= \int_0^{a^\dagger} r(a)\,(S(a,t)+I(a,t))\, da\, , \\
 W(t) &= \int_0^{a^\dagger} q(a)\, I(a,t)\, da\, .
 \end{aligned}
 $$
 Moreover, the following properties hold:
 \begin{itemize}
 \item[(i)] $\lim_{h\to 0} \frac{1}{h}[S(a+h,t+h)-S(a,t)]=
   -[K(a)W(t)+\mu(a)]S(a,t)$ a.e. in $[0,a^\dagger]\times\RR^+$,
 \item[(ii)] $\lim_{h\to 0} \frac{1}{h}[I(a+h,t+h)-I(a,t)]=
   K(a)W(t)S(a,t)-[\mu(a)+\alpha]I(a,t)$ a.e. in $[0,a^\dagger]
   \times\RR^+$,
 \item[(iii)] $\| S(\cdot,t)\|_{L^1(0,a^\dagger)}\leq(1+\alpha 
   a^\dagger)e^{R_0^d\|\beta\|_{L^\infty(0,a^\dagger)}(1+\alpha 
   a^\dagger)T}\| S_0+I_0\|_{L^1(0,a^\dagger)}$,
 \item[(iv)] $\| I(\cdot,t)\|_{L^1(0,a^\dagger)}\leq(1+\alpha 
   a^\dagger)e^{R_0^d\|\beta\|_{L^\infty(0,a^\dagger)}(1+\alpha 
   a^\dagger)T}\| S_0+I_0\|_{L^1(0,a^\dagger)}$,
 \item[(v)] there exist two constants $C_1,C_2>0$ depending on 
   $M$ and $T$ such that
   $$
   \| k(t)-\tilde k(t)\|_{(L^1(0,a^\dagger))^2} \leq C_1 e^{C_2t}
     \|(S_0,I_0)-(\tilde S_0,\tilde I_0)\|_{(L^1(0,a^\dagger))^2}
   $$
   where $\tilde k(t)$ is the solution relative to the initial 
   datum $(\tilde S_0,\tilde I_0)$.
 \end{itemize}
\end{theorem}

\section{Search for endemic equilibria}
\label{sec:endeq}

We now consider the problem of existence of steady states for the 
system (\ref{eq}).   Consequently, we are concerned with the problem
\begin{equation*}
 \left\{ \begin{aligned}
     &\frac{dS^*}{da}=-\left(\lambda^*(a)+\mu(a)\right)S^*(a)\, , \\
     &\frac{dI^*}{da}=\lambda^*(a)\,S^*(a)-\left(\mu(a)
       +\alpha\right)I^*(a)\, , \\
     &S^*(0)=R_0^d\Phi(Q^*)\int_0^{a^\dagger}\beta(a)\left(S^*(a)+I^*(a)
       \right)da\, , \\
     &I^*(0)=0\, ,
 \end{aligned}\right.
\end{equation*}
where
\begin{eqnarray}
 \label{lambda*}
 \lambda^*(a) & = &K(a)\int_0^{a^\dagger}q(\sigma)\, I^*(\sigma)\, 
   d\sigma\, ,\nonumber \\
 \label{Q*}
 Q^* & = & \int_0^{a^\dagger}r(\sigma)\left(S^*(\sigma)+I^*(\sigma)
   \right)d\sigma\, .
\end{eqnarray}

It is easy to see that
\begin{equation}
 \label{dfe}
 S^*_\mathrm{dfe}(a)=N^*(a)\, ,\qquad I^*_\mathrm{dfe}(a)\equiv 0
\end{equation}
is the disease-free equilibrium, provided that $N^*(a)$ is given by
(\ref{N*}).
Then we concentrate on the search of endemic states, that is nonnegative 
solutions for which $I^*(a)$ does not vanish identically.   Our aim in 
this Section is to reduce the solution of this problem to the solution 
of a system of equations of the scalar variables
\begin{eqnarray}
 \label{B}
 B^* &=& S^*(0)\, ,\\
 \label{W}
 W^* &=& \int_0^{a^\dagger}q(\sigma)\, I^*(\sigma)\, d\sigma\, .
\end{eqnarray}
To simplify, we denote
\begin{equation}
 \label{L}
 L(a)=\int_0^aK(\sigma)\, d\sigma\, 
\end{equation}
and, integrating, we obtain
\begin{equation}
 \label{S*,I*}
 \begin{aligned}
   S^*(a)&=B^*e^{-WL(a)}\pi(a)\, ,\\
   I^*(a)&=B^*W^*\pi(a)\int_0^aK(\sigma)e^{-W^*L(\sigma)-\alpha
     (a-\sigma)}\, d\sigma\, .
 \end{aligned}
\end{equation}

We define the following functions:
\begin{eqnarray*}
 F(\alpha,W) &=& \int_0^{a^\dagger}\beta(\sigma)\pi(\sigma)\left(
   e^{-\alpha\sigma}+\alpha\int_0^\sigma e^{-WL(\rho)-\alpha(\sigma-
   \rho)}d\rho\right)d\sigma\, ,\\
 G(\alpha,W) &=& \int_0^{a^\dagger}r(\sigma)\pi(\sigma)\left(
   e^{-\alpha\sigma}+\alpha\int_0^\sigma e^{-WL(\rho)-\alpha(\sigma-
   \rho)}d\rho\right)d\sigma\, ,\\
\ H(\alpha,W) &=& \int_0^{a^\dagger}q(\sigma)\pi(\sigma)\int_0^\sigma 
   K(\rho)e^{-WL(\rho)-\alpha(\sigma-\rho)}d\rho\, d\sigma\, .
\end{eqnarray*}
We note that, integrating by parts, $F(\alpha,W)$ and $G(\alpha,W)$ can 
be written as
\begin{equation}
 \label{F,G}
 \begin{aligned}
   F(\alpha,W) &= \int_0^{a^\dagger}\beta(\sigma)\pi(\sigma)\left(
     e^{-WL(\sigma)}+W\int_0^\sigma K(\rho)e^{-WL(\rho)-\alpha(\sigma-
     \rho)}d\rho\right)d\sigma\, ,\\
   G(\alpha,W) &= \int_0^{a^\dagger}r(\sigma)\pi(\sigma)\left(
     e^{-WL(\sigma)}+W\int_0^\sigma K(\rho)e^{-WL(\rho)-\alpha(\sigma-
     \rho)}d\rho\right)d\sigma\, .
 \end{aligned}
\end{equation}

Substituting (\ref{S*,I*}) in (\ref{W}) we get
\begin{equation}
 \label{W2}
 \begin{aligned}
 W^* &=B^*W^*\int_0^{a^\dagger}q(a)\pi(a)\int_0^aK(\sigma)
      e^{-W^*L(\sigma)-\alpha(a-\sigma)}d\sigma\, da\\
     &=B^*W^*H(\alpha,W^*)\, .
 \end{aligned}
\end{equation}

Substituting (\ref{S*,I*}) in (\ref{Q*}), we have
\[
Q^* = B^*\int_0^{a^\dagger}r(a)\pi(a)\left(e^{-W^*L(a)}+W^*\int_0^a
 K(\sigma)e^{-W^*L(\sigma)-\alpha(a-\sigma)}d\sigma\right)da
\]
and, using (\ref{F,G}) and (\ref{W2}),
\begin{equation}\label{qu*}
Q^* = \frac{G(\alpha,W^*)}{H(\alpha,W^*)}\, .
\end{equation}
Now our problem has been reduced to solving the following system of 
equations:
\begin{equation}
 \label{sys}
 \left\{\begin{aligned}
     &B^*=\frac{1}{H(\alpha,W^*)}\, ,\\
     &R_0^d\Phi\left(\frac{G(\alpha,W^*)}{H(\alpha,W^*)}\right)
       F(\alpha,W^*)=1\, .
   \end{aligned}\right.
\end{equation}
The second equation in (\ref{sys}) depends only on $W^*$:  
once we have a solution of it we can substitute it in the first equation.  
So, in order to simplify the second equation, we write
\begin{equation}
 \label{varphi}
 \varphi(\alpha,W)=R_0^d\Phi\left(\frac{G(\alpha,W)}{H(\alpha,W)}\right)
   F(\alpha,W)
\end{equation}
and we are left with equation
$$
\varphi(\alpha,W)=1\, .
$$

We have that
$$
\varphi(\alpha,0)=R_0^d\Phi\left(\frac{\int_0^{a^\dagger}r(\sigma)
   \pi(\sigma)\, d\sigma}{H(\alpha,0)}\right)
$$
and 
$$
\lim_{x\to +\infty}\varphi(\alpha,x)=0\, ,
$$
so the second equation in (\ref{sys}) has at least one solution 
$W^*$ if $\varphi(\alpha,0) >1$. Since $\Phi(\cdot)$ is decreasing in 
$[0,Q^*_d]$, this latter condition is equivalent to
$$
\frac{\int_0^{a^\dagger}r(\sigma)\pi(\sigma)\, d\sigma}{H(\alpha,0)} 
 <Q_d^*\, ,
$$
where $Q_d^*$ is the solution of (\ref{Qd}), that is to (see (\ref{N*}))
\begin{equation}\label{condizione}
\frac{Q_d^* H(\alpha,0)}{\int_0^{a^\dagger}r(\sigma)\pi(\sigma)
 \, d\sigma}=N^*(0) H(\alpha,0)>1\, .
\end{equation}
Condition (\ref{condizione}) can be written as 
\begin{equation}\label{mag1}
R_0^e > 1\, ,
\end{equation}
where $R_0^e$ is the basic epidemic reproduction ratio, i.e. the 
number of secondary cases which one case would produce in a 
completely susceptible population.    This number is defined in this 
context as the spectral radius of the next-generation operator $\mathcal{G}$ 
(see \cite{DH}, Chapter 7)
$$
(\mathcal{G}u)(a) = \int_0^{a^\dagger} g(a,\sigma)\, u(\sigma)\, d\sigma\, ,
$$
where
$$
g(a,\sigma) = \int_\sigma^{a^\dagger} K(a)q(\rho)N^*(a)
  \frac{\pi(\rho)}{\pi(\sigma)} e^{-\alpha(\rho-\sigma)}\, d\rho\, .
$$
It is easy to see that
$$
R_0^e = \int_0^{a^\dagger}  K(a)N^*(a) \int_a^{a^\dagger} q(\sigma)
  \frac{\pi(\sigma)}{\pi(a)} e^{-\alpha(\sigma-a)}\, d\sigma\, da\, .
$$

Since $N^*(0)=\frac{N^*(a)}{\pi(a)}$ for any $a\in[0,a^\dagger]$, we 
can write
$$
\begin{aligned}
N^*(0)H(\alpha,0) = &\int_0^{a^\dagger} q(\sigma)\pi(\sigma)\int_0^\sigma 
   K(a)e^{-\alpha(\sigma-a)}\frac{N^*(a)}{\pi(a)}da\, d\sigma \\
 = &\int_0^{a^\dagger}K(a)\frac{N^*(a)}{\pi(a)}\int_a^{a^\dagger} 
   q(\sigma)\pi(\sigma)e^{-\alpha(\sigma-a)}d\sigma\, da = R_0^e\, .
\end{aligned}
$$

Thus (\ref{mag1}) is a sufficient condition for the existence of an endemic 
equilibrium.

Next section is devoted to discuss uniqueness of equilibria under such 
condition.

\section{Discussing uniqueness of an endemic equili\-brium}\label{secexun}

Condition (\ref{mag1}) is sufficient to have at least one endemic state, 
but for uniqueness we need some additional assumptions.
A first uniqueness case occurs when the disease does not induce mortality 
(see also \cite{CIM}).

\begin{theorem}
  \label{trm-uniqueness1}
  If $\alpha=0$ and $R_0^e>1$, there exists one and only one endemic 
  equilibrium of the problem (\ref{eq}).
\end{theorem}

\begin{proof}
  Existence of endemic equilibria is equivalent to existence of 
  solutions of the system (\ref{sys}).   Since the function
  $$
  \frac{\int_0^{a^\dagger}r(\sigma)\pi(\sigma)
    \, d\sigma}{\int_0^{a^\dagger}q(\sigma)\pi(\sigma)\int_0^\sigma 
    K(\rho)e^{-WL(\rho)}d\rho\, d\sigma}
  $$
  is increasing in $W$ and tends to infinity as $W$ tends to infinity, 
  there exists $\widetilde W$ such that
  $$
  \varphi(0,W)=\Phi\left(\frac{\int_0^{a^\dagger}r(\sigma)\pi(\sigma)
    \, d\sigma}{\int_0^{a^\dagger}q(\sigma)\pi(\sigma)\int_0^\sigma 
    K(\rho)e^{-WL(\rho)}d\rho\, d\sigma}\right)\int_0^{a^\dagger}
    \beta(\sigma)\pi(\sigma)\, d\sigma
  $$
  is decreasing for $W$ in $[0,\widetilde W]$, with $\varphi(0,
  \widetilde W)<1$, and non-increasing for $W>\widetilde W$.   Now, 
  $R_0^e>1$ implies $\varphi(0,0)>1$ and so there exists one and only 
  one solution of $\varphi(0,W)=1$.
\end{proof}

The proof of another condition for $\alpha > 0$, requires the following 
Lemma (see \cite{HadelerDietz})

\begin{lemma}[Hadeler and Dietz]
 \label{lmm-HadelerDietz}
 Let $g$, $u$ and $v$ be locally summable real functions on the (finite 
 or infinite) interval $(a,b)$, with $u(x),\, v(x)\geq 0$, such that
 \begin{equation}
   \label{hypHD}
   x\leq y \Longrightarrow g(x)\leq g(y)\quad\mbox{ and }\quad v(x)u(y)
     \geq u(x)v(y)\, .
 \end{equation}
 Then
 \begin{equation}
   \label{claimHD}
   \int_a^b g(x)u(x)\, dx \cdot \int_a^b v(x)\, dx \geq 
     \int_a^b g(x)v(x)\, dx \cdot \int_a^b u(x)\, dx\, ,
 \end{equation}
 provided the integrals exist.

 Moreover, when $u$, $v$ and $g$ are positive and $\frac{u(x)}{v(x)}$ 
 is strictly increasing, then (\ref{claimHD}) holds with strict 
 inequality.
\end{lemma}

Then we obtain the following result.

\begin{theorem}
 \label{trm-uniqueness2}
 If for any $0\leq\rho_1\leq\rho_2\leq a^\dagger$,
 \begin{equation}
   \label{hyp}
   \begin{aligned}
     K(\rho_1) &\int_{\rho_1}^{a^\dagger}q(\sigma)\pi(\sigma)
        e^{-\alpha\sigma}d\sigma \int_{\rho_2}^{a^\dagger}r(\sigma)
        \pi(\sigma)e^{-\alpha\sigma}d\sigma \\
       &\leq K(\rho_2)\int_{\rho_1}^{a^\dagger}r(\sigma)\pi(\sigma)
        e^{-\alpha\sigma}d\sigma \int_{\rho_2}^{a^\dagger}q(\sigma)
        \pi(\sigma)e^{-\alpha\sigma} d\sigma\, .
   \end{aligned}
 \end{equation}
 and $R_0^e>1$, then there exists one and only one endemic equilibrium 
 of the problem (\ref{eq}).
\end{theorem}
\begin{proof}
  Fixed $\alpha\geq 0$, $R_0^e>1$ implies $\varphi(\alpha,0)>1$.  By 
  definition $F(\alpha,W)$ is decreasing in $W$.  Then, if we prove 
  that $\frac{G(\alpha,W)}{H(\alpha,W)}$ is increasing in $W$, since 
  $\lim_{W\to +\infty}\frac{G(\alpha,W)}{H(\alpha,W)}=+\infty$, we have 
  that $\varphi(\alpha,\cdot)$ is decreasing on a suitable interval 
  $[0,\widetilde W]$, where $\phi\left(\frac{G(\alpha,\widetilde W)}
  {H(\alpha,\widetilde W)}\right)>0$, and this gives our claim.   Since
  $$
  \frac{\partial}{\partial W}\left(\frac{G(\alpha,W)}{H(\alpha,W)}\right)
    =\frac{\frac{\partial G(\alpha,W)}{\partial W}H(\alpha,W)-G(\alpha,W)
    \frac{\partial H(\alpha,W)}{\partial W}}{(H(\alpha,W))^2}\, ,
  $$
  we study the sign of the numerator $G_WH-GH_W$.  If we denote
  $$
  u(\rho)=K(\rho)e^{-WL(\rho)+\alpha\rho}\int_\rho^{a^\dagger}q(\sigma)
    \pi(\sigma)e^{-\alpha\sigma} d\sigma
  $$
  and
  $$
  v(\rho)=e^{-WL(\rho)+\alpha\rho}\int_\rho^{a^\dagger}r(\sigma)
    \pi(\sigma)e^{-\alpha\sigma}d\sigma\, ,
  $$
  we get:
  $$
  \begin{aligned}
  G_WH-GH_W = &-\alpha\int_0^{a^\dagger} L(\rho)v(\rho)\, d\rho\,
       \int_0^{a^\dagger}u(\rho)\, d\rho\\
     &+\alpha\int_0^{a^\dagger}v(\rho)\, d\rho\,\int_0^{a^\dagger} 
       L(\rho)u(\rho)\, d\rho\\
     &+\int_0^{a^\dagger}r(\sigma)\pi(\sigma)e^{-\alpha\sigma} 
      d\sigma\, \int_0^{a^\dagger} L(\rho)u(\rho)\, d\rho\, .
  \end{aligned}
  $$
  Now we apply Hadeler-Dietz's lemma.   Since $L(\rho)$ by definition 
  (see (\ref{L})) is increasing, we choose $g\equiv L$.   The 
  second part of condition (\ref{hypHD}) coincides with 
  (\ref{hyp}), so we obtain that
  $$
  G_WH-GH_W\geq\int_0^{a^\dagger}r(\sigma)\pi(\sigma)e^{-\alpha\sigma} 
    d\sigma\, \int_0^{a^\dagger} L(\rho)u(\rho)\, d\rho > 0.
  $$
\end{proof}

\begin{remark}
  \label{rmk-uni}
  Condition (\ref{hyp}) is verified for example when there exists 
  a function $h(\cdot)$ such that $q(a)= h(a)r(a)$ and both $K(\cdot)$ 
  and $h(\cdot)$ are non-decreasing functions.
\end{remark}

We analyzed so far sufficient conditions to have uniqueness. We can also 
show examples in which multiple equilibria occur. The following case 
corresponds to a very special situation but it helps to understand the 
mechanisms responsible for non-uniqueness.

\begin{figure}[h]
 \begin{center}
   \includegraphics[width=6cm]{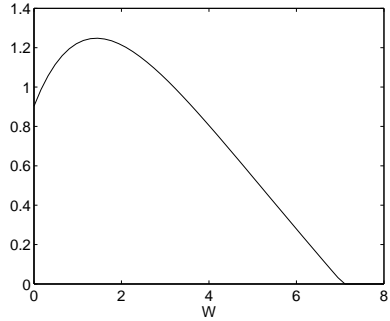}
   \caption{the graph of $\varphi(10,W)$ for the choices (\ref{choices}).}
   \label{fig1}
 \end{center}
\end{figure}

Let us consider the following choices
\begin{equation}
 \label{choices}
 \begin{array}{l}
   a^\dagger=\frac{\pi}{2} \\
   \alpha=10 \\
   R_0^d=27 \\
   \beta(a)\equiv 1 \\
   q(a)\equiv 1 \\
 \end{array}
 \hspace{2cm}
 \begin{array}{l}
   r(a)\equiv 1 \\
   \mu(a)=\tan a \\
   K(a)=\left\{\begin{array}{ll}
       1 & \mbox{if }a\in\left[0,\frac{\pi}{6}\right]\cup\left[\frac{\pi}
         {3},\frac{\pi}{2}\right] \\
       0 & \mbox{if }a\in\left(\frac{\pi}{6},\frac{\pi}{3}\right)
     \end{array}\right. \\
   \Phi(x)=\max\left\{1-\frac{x}{18},0\right\}
 \end{array}
\end{equation}
Then we have
$$
\pi(a)=\cos(a)\, ,\qquad 
L(a)=\left\{\begin{array}{ll}
     a & \mbox{if }a\in\left[0,\frac{\pi}{6}\right], \\
     \frac{\pi}{6} & \mbox{if }a\in\left(\frac{\pi}{6},
       \frac{\pi}{3}\right), \\
     a-\frac{\pi}{6} & \mbox{if }a\in\left[\frac{\pi}
       {3},\frac{\pi}{2}\right].
   \end{array}\right.
$$
Thus we may explicitly compute the functions $F$, $G$ and $H$ and we 
find that the function (see (\ref{varphi}))
\[
\varphi(10,W)=R_0^d\Phi\left(\frac{G(10,W)}{H(10,W)}\right)F(10,W)
\]
has the graph shown in Figure \ref{fig1}. This leads to the following result

\begin{theorem}
  Problem (\ref{eq}) with the choices (\ref{choices}) has two endemic 
  equilibria.
\end{theorem}

\begin{figure}[h]
 \begin{center}
   \includegraphics[width=6cm]{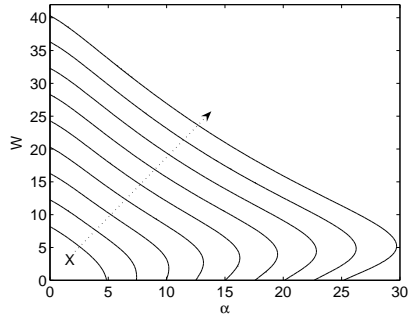}
   \caption{endemic equilibria curves 
     for $X=10,14,18,22,26,30,34,38,42$ in (\ref{plus}).}
   \label{fig2}
 \end{center}
\end{figure}

As an extension of this special case we consider the following forms
\begin{equation}\label{plus}
R_0^d=\frac{3X}{2}\quad\mbox{and}\quad\Phi(x)=\max\left\{1-
 \frac{x}{X},0\right\}
\end{equation}
and in Figure \ref{fig2} we show the endemic equilibria curves versus the 
parameter $\alpha$, for different choices of $X$.
The special case of (\ref{choices}) corresponds to $X=18$. The 
behaviour becomes more evident as  $X$ increases.
%

\section{Characteristic equation}\label{sec:ce}

In order to investigate the local asymptotic stability of the 
equilibria we consider the characteristic equation.  For this standard 
technique see for example \cite{Iannelli}.  In the following sections we 
assume that $\Phi$ is differentiable at $Q^*_d$ (in the case of the 
disease free equilibrium) and in $Q^*$ defined in (\ref{Q*}) (in the 
case of the endemic equilibria).

A convenient way to represent problem (\ref{eq}) in an equivalent form 
is to use the variables $B(t)$ as in (\ref{B(t)}) and $W(t)$ as in 
(\ref{W(t)}). In fact, using (\ref{Scar}) and (\ref{Icar}), we get 
a system of two integral equations that, for $t\geq a^\dagger$, is given by
\begin{equation}
 \label{BW}
 \left\{\hspace{-1mm}\begin{aligned}
   B&(t)=R_0^d\Phi\left(\int_0^{a^\dagger}r(a)\pi(a)\left(e^{-\int_0^a
       K(\sigma)W(t-a+\sigma)\, d\sigma}\phantom{\int_0^a}\right.\right.\\
    &\hspace{-2mm}\left.\left.+e^{-\alpha a} \hspace{-2mm} \int_0^a \hspace{-2mm}K(\sigma)e^{\alpha\sigma-
       \int_0^\sigma K(\rho)W(t-a+\rho)\, d\rho}W(t-a+\sigma)\, 
       d\sigma\right)\hspace{-1mm}B(t-a)\,da\right)\hspace{-1mm}\cdot\\
    &\hspace{-2mm}\cdot\int_0^{a^\dagger}\beta(a)\pi(a)\left(e^{-\int_0^a
       K(\sigma)W(t-a+\sigma)\, d\sigma}\phantom{\int_0^a}\right.\\
    &\hspace{-2mm}\left.+e^{-\alpha a}\int_0^a \hspace{-2mm}K(\sigma)e^{\alpha\sigma-
       \int_0^\sigma K(\rho)W(t-a+\rho)\, d\rho}W(t-a+\sigma)\, 
       d\sigma\right) \hspace{-1mm}B(t-a)\,da\, ,\\
   W&(t)=\hspace{-1mm}\int_0^{a^\dagger} \hspace{-3mm}q(a)\pi(a)e^{-\alpha a}
       \hspace{-1mm}\int_0^a 
       \hspace{-3mm}K(\sigma)e^{\alpha\sigma-\int_0^\sigma K(\rho)W(t-a+\rho)\, 
       d\rho}W(t-a+\sigma)\, d\sigma\cdot\\
    &\hspace{-2mm}\cdot B(t-a)\,da\, .
   \end{aligned}\right.
\end{equation}
Actually, the limiting equations of the two general integral equations for 
$B(t)$, $W(t)$ coincide with (\ref{BW}).  Moreover, the constant solutions 
$B(t)\equiv B^*$, $W(t)\equiv W^*$ of (\ref{BW}) correspond to the 
different steady states of the problem. Namely we have:
\begin{itemize}
\item the trivial solution $B^*=W^*=0$;
\item the disease-free equilibrium corresponding to $W^*=0$ and $B^*$ 
  provided by the equation
  \begin{equation}\label{sys-dfe}
  1=R_0^d\Phi\left(B^*\int_0^{a^\dagger}r(a)\pi(a)\,da\right);
  \end{equation}
\item the endemic equilibria provided by system (\ref{sys}).
\end{itemize}
Note that the case of the disease-free equilibrium corresponds to 
(\ref{dfe}) with (\ref{N*}) and (compare (\ref{sys-dfe}) with  (\ref{Qd}))
\begin{equation}
\label{QdB*}
Q_d^*= B^* \int_0^{a^\dagger}r(a)\pi(a)\,da.
\end{equation}
Finally, the endemic equilibrium is given by (\ref{S*,I*}).

In order to linearise (\ref{BW}) at $(B^*,W^*)$, 
we let
$$
b(t)=B(t)-B^* \quad\mbox{ and }\quad w(t)=W(t)-W^*\, .
$$
Then we obtain the linear system
\begin{equation}
 \label{bw}
 \left\{\begin{aligned}
   b(t) &= \int_0^{a^\dagger}\Psi_1(a)b(t-a)\, da+ 
         \int_0^{a^\dagger}\Psi_2(a)w(t-a)\, da\, ,\\
   w(t) &= \int_0^{a^\dagger}\Psi_3(a)b(t-a)\, da+ 
         \int_0^{a^\dagger}\Psi_4(a)w(t-a)\, da\, ,
 \end{aligned}\right.
\end{equation}
where the convolution kernels $\Psi_1$, $\Psi_2$, $\Psi_3$ and $\Psi_4$ 
are given by
\begin{equation}
\begin{aligned}
 \Psi_1(a)= & \left( R_0^d\Phi(Q^*)\beta(a) + B^*
         \frac{\Phi'(Q^*)}{\Phi(Q^*)}r(a) \right)\cdot \\
        &\cdot  \pi(a) \left( e^{-W^*L(a)} 
     + W^*\int_0^aK(\sigma)e^{-\alpha(a-\sigma)
         -W^*L(\sigma)} d\sigma \right), \\
 \Psi_2(a)= &-\alpha B^* \int_a^{a^\dagger}  \left( R_0^d
         \Phi(Q^*)\beta(\sigma) + B^*\frac{\Phi'(Q^*)}{\Phi(Q^*)}
         r(\sigma) \right) \pi(\sigma)K(\sigma-a) \cdot \\
       &\cdot \int_{\sigma-a}^\sigma e^{-\alpha(\sigma-\rho)-
         W^*L(\rho)}d\rho\, d\sigma\, , \\
 \Psi_3(a)= &W^*q(a)\pi(a)\int_0^aK(\sigma)e^{-\alpha(a-\sigma)
         -W^*L(\sigma)}d\sigma\, , \\
 \Psi_4(a)= &B^*\int_a^{a^\dagger}q(\sigma)\pi(\sigma)K(\sigma-a)\cdot\\
         & \cdot\left( e^{-\alpha a-W^*L(\sigma-a)} 
       - W^*\int_{\sigma-a}^\sigma K(\rho)e^{-\alpha(\sigma-\rho)
         -W^*L(\rho)} d\rho \right) d\sigma\, ,
\end{aligned}
\label{convker}
\end{equation}
with $Q^*$ defined in (\ref{Q*}).

Taking Laplace transforms in (\ref{bw}), we obtain the following 
characteristic equation (see \cite{Iannelli}):
\begin{equation}
 \label{eqcar}
 \Psi(\lambda)=(1-\widehat\Psi_1(\lambda))(1-\widehat\Psi_4(\lambda))-
   \widehat\Psi_2(\lambda)\widehat\Psi_3(\lambda)=0\, ,
\end{equation}
where for $i=1,2,3,4$
$$
\widehat\Psi_i(\lambda)=\int_0^{+\infty} e^{-\lambda a}\Psi_i(a)\, da
$$
represents the Laplace transform of $\Psi_i(a)$.

To study the stability of an equilibrium,   we need to determine the 
location of the roots of (\ref{eqcar}). In the following section we will 
be concerned with the disease-free equilibrium, then we will investigate 
the role of the parameter $\alpha$.

\section{Stability of the disease-free equilibrium}\label{sec:stabdfe}

In the case of the disease-free equilibrium (\ref{dfe}), the convolution 
kernels become
$$
\begin{aligned}
 \Psi_1(a)= &\beta(a)\pi(a)+R_0^dB^*\Phi'(Q^*_d)\, r(a)\pi(a)\, ,\\
 \Psi_2(a)= &B^*(e^{-\alpha a}-1)\int_a^{a^\dagger} 
          \beta(\sigma)\pi(\sigma)K(\sigma-a)\, d\sigma\\
         &+R_0^d(B^*)^2\Phi'(Q^*_d)\int_0^{a^\dagger}\hspace{-0.3cm}
          \beta(\sigma)\pi(\sigma)\, d\sigma(e^{-\alpha a}-1)
          \int_a^{a^\dagger}\hspace{-0.3cm} r(\sigma)\pi(\sigma)
          K(\sigma-a)\, d\sigma\, ,\\
 \Psi_3(a)= &0\, ,\\
 \Psi_4(a)= &B^*e^{-\alpha a}\int_a^{a^\dagger}q(\sigma)\pi(\sigma)
          K(\sigma-a)\, d\sigma\, ,
\end{aligned}
$$
where $Q^*_d$ is given in (\ref{QdB*}).

Hence, the roots of (\ref{eqcar}) are the union of the roots of the two 
equations
\[
\hat\Psi_1(\lambda)=1,\quad \hat\Psi_4(\lambda)=1,
\]
that can be considered separately. Note that the first of this 
equations is exactly the characteristic equation (\ref{carN}) 
for the demographic problem (\ref{Gurtin_MC}) in the absence of 
the disease. Thus we have

\begin{theorem}\label{trm-stabdfe}
  If $R_0^e>1$, the disease-free equilibrium is unstable.
  If $R_0^e<1$, the disease-free equilibrium is 
  stable or unstable depending on whether the equilibrium (\ref{N*}) of the 
  population is stable or not for (\ref{Gurtin_MC}).
\end{theorem}

\begin{proof}
  Let $R_0^e>1$.  We have that
  $$
  \begin{aligned}
  \hat\Psi_4(0)&=\int_0^{a^\dagger} \Psi_4(a)\, da=B^*\int_0^{a^\dagger}
      e^{-\alpha a}\int_a^{a^\dagger}q(\sigma)\pi(\sigma)K(\sigma-a)
      \, d\sigma\, da\\
    &=B^*H(\alpha,0)\, .
  \end{aligned}
  $$
  Since
  $$
  \frac{G(\alpha,0)}{H(\alpha,0)} = \frac{Q^*_d}{N^*(0)H(\alpha,0)} < Q^*_d
  $$
  and $\Phi(\cdot)$ is decreasing on $[0,Q^*_d]$, one has
  $$
  \varphi(\alpha,0)=R_0^d\Phi\left(\frac{G(\alpha,0)}{H(\alpha,0)}
    \right)>1\, .
  $$
  Now, comparing (\ref{sys-dfe}) with the previous equation, we get that
  $$
  B^*>\frac{1}{H(\alpha,0)}
  $$
  and consequently $\hat\Psi_4(0)>1$.   This means that $1-
  \hat\Psi_4(t)$ has exactly one root on the positive real line.
  
  If $R_0^e<1$, then $\hat\Psi_4(0)<1$ and, if $\lambda\in\CC$ has 
  positive real part,
  $$
  \begin{aligned}
  |1-\hat\Psi_4(\lambda)| &\geq 1 - \left|\int_0^{a^\dagger}
      e^{-\lambda a}\Psi_4(a)\, da\right|\\
    &> 1-\int_0^{a^\dagger} \Psi_4(a)\, da= 1-\hat\Psi_4(0)>0\, .
  \end{aligned}
  $$
  So, there are no zeros of $1-\hat\Psi_4(\lambda)$ with positive real 
  part.   The study of the zeros of $1-\hat\Psi_1(\lambda)$ coincides 
  with the stability analysis of the non trivial equilibrium of the 
  total population (equation (\ref{Gurtin_MC})) and this 
  concludes the proof. 
\end{proof}

From the previous Theorem we see that the epidemic reproduction ratio 
$R_0^e$ determines also stability of the disease free equilibrium.
In the next Section we discuss how instability may depend on the 
disease induced mortality $\alpha$.

\section{Destabilizing effect of the extra-mortality}\label{secsette}

Since an analytic study of the characteristic equation at the endemic 
equilibrium does not seem to be possible in the general case, we 
present a specific example. We consider a special set of parameters 
that shows how the extra-mortality $\alpha$ can destabilize the endemic 
equilibrium.

Let us first consider the following choices:
\begin{equation}
 \label{choices-stab}
   a^\dagger=\frac{\pi}{2}\, , \qquad
   \beta(a)=\frac{3}{2}\sin(2a)\, , \qquad
   r(a)=\frac{3}{2}\sin(2a)\, , \qquad
   \mu(a)=\tan a\, . 
\end{equation}

With these choices and $\alpha=0$ the convolution kernels (\ref{convker}) 
for the endemic equilibrium  become
$$
\begin{aligned}
 \Psi_1(a)= & \left( 1+R_0^dB^*\Phi'(Q^*) \right) \frac{3}{2}\sin(2a)
   \cos(a)\, ,\\
 \Psi_2(a)= & 0\, , \\
 \Psi_3(a)= &q(a)\cos(a) \left( 1-e^{-W^*L(a)} \right)\, , \\
 \Psi_4(a)= &B^*\int_a^{a^\dagger}q(\sigma)\cos(\sigma)K(\sigma-a)
         e^{-W^*L(\sigma)} d\sigma\, ,
\end{aligned}
$$
where $B^*$ and $Q^*$ are given by (\ref{sys}) and (\ref{qu*}).
Then we have the following result.

\begin{proposition}
\label{prp-radim}
  Let (\ref{choices-stab}) be satisfied.   If $\alpha=0$ and
  \begin{equation}
    \label{hyp-radim}
    R_0^dB^*\Phi'(Q^*) = -5\, ,
  \end{equation}
  then the characteristic equation (\ref{eqcar}) has two  imaginary 
  roots  $\lambda_\pm =\pm 5i$  and any other root has negative real part.
\end{proposition}

\begin{proof}
 Since $\Psi_2(a)=0$, the characteristic equation is simply
 $$
 \Psi(\lambda)=(1-\widehat\Psi_1(\lambda))(1-\widehat\Psi_4(\lambda))=0\, .
 $$
 If $\lambda$ has nonnegative real part, we have
 \begin{eqnarray*}
   \left|\Re\left(\widehat\Psi_4(\lambda) \right)\right| 
     & \leq & \left| B^*\int_0^{a^\dagger}e^{-\lambda a}
       \int_a^{a^\dagger}q(\sigma)\pi(\sigma)K(\sigma-a)
       e^{-W^*L(\sigma)} d\sigma\, da \right| \\
     & \leq &  B^* \int_0^{a^\dagger} \int_a^{a^\dagger}q(\sigma)
       \pi(\sigma)K(\sigma-a) e^{-W^*L(\sigma)} d\sigma\, da\, .
 \end{eqnarray*}
 Since
 $$
 \begin{aligned}
   \int_0^{a^\dagger} \hspace{-0.2cm} &\int_a^{a^\dagger}
       \hspace{-0.2cm} q(\sigma)\pi(\sigma)K(\sigma-a)
       e^{-W^*L(\sigma)} d\sigma\, da = \int_0^{a^\dagger} 
       \hspace{-0.2cm}q(\sigma)\pi(\sigma)e^{-W^*L(\sigma)} 
       \int_0^\sigma \hspace{-0.2cm} K(\rho)\, d\rho\,d\sigma \\
     &< \int_0^{a^\dagger} \hspace{-0.2cm}q(\sigma)\pi(\sigma) 
       \int_0^\sigma \hspace{-0.2cm} K(\rho) e^{-W^*L(\rho)}\, 
       d\rho\,d\sigma = H(0,W^*) = \frac{1}{B^*}\, ,
 \end{aligned}
 $$
 $\left|\Re\left(\widehat\Psi_4(\lambda) \right)\right|<1$.   
 Then the roots of the characteristic equation with nonnegative 
 real part can be found only if $\widehat\Psi_1(\lambda)=1$, i.e.
 \begin{equation}
 \label{Psi1=1}
 (1+\tau)\frac{3}{2}\int_0^\frac{\pi}{2}e^{-\lambda a}\sin(2a)
   \cos(a)\, da = 1
 \end{equation}
 with $\tau=R_0^dB^*\Phi'(Q^*)$.   We look for values of 
 $\tau<0$ for which a couple of imaginary roots $\lambda
 =\pm i\omega$ of (\ref{Psi1=1}) exist.  This is equivalent to 
 solving the following system:
 $$
 \left\{ \begin{array}{l}
       \int_0^\frac{\pi}{2} \sin(\omega a)\sin(2a)\cos(a)\, da=0\, ,\\
       (1+\tau)\frac{3}{2}\int_0^\frac{\pi}{2} \cos(\omega a)
         \sin(2a)\cos(a)\, da=1\, .
       \end{array} \right.
 $$
 The first equation gives
 $$
 \frac{\sin\frac{\pi(\omega-3)}{2}}{4(\omega-3)}+
   \frac{\sin\frac{\pi(\omega-1)}{2}}{4(\omega-1)}-\frac{\sin
   \frac{\pi(\omega+1)}{2}}{4(\omega+1)}-\frac{\sin\frac{\pi
   (\omega+3)}{2}}{4(\omega+3)} = \frac{4\omega\cos
   \frac{\pi\omega}{2}}{(\omega^2-9)(\omega^2-1)} = 0
 $$
 and the solutions are $\omega_k^\pm=\pm(2k+1)$ with 
 $k=2,3,\dots$.   Correspondingly, the values of $\tau=
 \tau_k$ are given by the following equality
 $$
 \begin{aligned}
 \tau_k &= \frac{2}{3\int_0^\frac{\pi}{2}\cos((2k+1)a)\sin(2a)
     \cos(a)\, da}-1 \\
   &= \frac{2}{\frac{3}{8} \left[ \frac{\cos(2(k-1)a)}{k-1} + 
     \frac{\cos(2ka)}{k} - \frac{\cos(2(k+1)a)}{k+1} - 
     \frac{\cos(2(k+2)a)}{k+2} \right]_0^\frac{\pi}{2}}-1\, , \\
 \end{aligned}
 $$
 which means
 $$
 \tau_k = \left\{ \begin{array}{ll}
     \frac{1-4k^2}{3} & \mbox{for $k$ even}\, , \\
     \frac{-4k^2-8k-3}{3} & \mbox{for $k$ odd}\, .
   \end{array} \right.
 $$
 It is easy to verify that $\tau_k$ is decreasing with $k$.   
 Then $\tau_2=-5$ is the maximum of the $\tau_k$'s.   
 Since by the implicit function theorem one obtains that 
 $\frac{\partial\Re\lambda}{\partial\tau}(-5)<0$, 
 $\lambda=\omega_2^\pm i=\pm 5i$ are the first roots 
 to cross the imaginary axis and this completes the proof.
\end{proof}

We now let  $\alpha$ become positive, and analyze how  the stability 
of the endemic equilibrium changes with $\alpha$. For simplicity we 
make a further choice considering 
\begin{equation}
\label{choices-stab2}
\Phi(x)=\max\left\{1-\frac{x}{10},0\right\}\, , \qquad q(a)=\frac{3}{2}
  \sin(2a)\, , \qquad K(a)=a\, .
\end{equation}
By (\ref{choices-stab}) and Remark \ref{rmk-uni}, for every fixed 
$\alpha$ there exists one and only one endemic equilibrium.   So, 
equation $\varphi(\alpha,W)=1$ has only one solution $(0,W^*_0)$, 
with $\alpha=0$.   With the choices made in (\ref{choices-stab}), 
(\ref{hyp-radim}) and (\ref{choices-stab2}), we get
$$
R_0^d=6 \qquad \mbox{and} \qquad H(0,W^*_0)=\frac{3}{25}\, .
$$
so that  $W^*_0\simeq 5.04512$. Moreover at this point we have 
$\displaystyle\frac{\partial\varphi(0,W^*_0)}{\partial W}\neq 0$.   
Then for $\alpha>0$ sufficiently small, there exists a function $
W^*(\alpha)$ such that 
\begin{equation}\label{dini}
\varphi(\alpha,W^*(\alpha))=1
\end{equation}
and we obtain a branch of endemic states $\left(B^*(\alpha),
W^*(\alpha)\right)$.

Now we show that this endemic equilibrium changes its stability, 
in the sense that at $\alpha=0$ the roots of the 
characteristic equation cross forward the imaginary axis.

We write $\lambda=\zeta+i\omega$, so that we can consider the 
characteristic equation (\ref{eqcar}) at the equilibrium $\left(B^*(\alpha),W^*(\alpha)\right)$, as a system depending on $\alpha$
$$
\begin{aligned}
 F_1(\alpha,\zeta,\omega) &= \Re\Psi(\lambda) = 0\, , \\
 F_2(\alpha,\zeta,\omega) &= \Im\Psi(\lambda) = 0\, .
\end{aligned}
$$

By Proposition~\ref{prp-radim} we know that $F_1(0,0,\pm 5)=
F_2(0,0,\pm 5)=0$, while detailed calculations show 
$$
\begin{aligned}
\frac{\partial F_1}{\partial\omega}(0,0,\pm 5) &= - 
 \frac{\partial F_2}{\partial\zeta}(0,0,\pm 5), \\
\frac{\partial F_2}{\partial\omega}(0,0,\pm 5) &= 
 \frac{\partial F_1}{\partial\zeta}(0,0,\pm 5),
\end{aligned}
$$
so that the jacobian of the system with respect to $\zeta$ and $\omega$ 
at $(0,0,\pm 5)$ is equal to
\[
\left(\frac{\partial {F_1}}{\partial\zeta}(0,0,\pm 5)\right)^2+\left(
  \frac{\partial F_2}{\partial\zeta}(0,0,\pm 5)\right)^2
\]
and by numerical computation (the evaluations are obtained using 
Mathematica, \texttt{www.wolfram.com}) we conclude that it is positive.
Then for $\alpha>0$ sufficiently 
small there exist $\zeta(\alpha)$ and $\omega(\alpha)$ such that 
$\zeta(0)=0$, $\omega(0)=\pm 5$ and $F_i(\alpha,\zeta(\alpha),
\omega(\alpha))=0$ for $i=1,2$. Concerning the sign 
of $\zeta'(0)$,  again by numerical computation we have
$$
\zeta'(0) = \frac{-\frac{\partial F_1}{\partial\alpha}(0,0,\pm 5)
 \frac{\partial F_1}{\partial\zeta}(0,0,\pm 5)-\frac{\partial 
 F_2}{\partial\alpha}(0,0,\pm 5)\frac{\partial F_2}{\partial\zeta}
 (0,0,\pm 5)}{\left(\frac{\partial F_1}{\partial\zeta}(0,0,\pm 5)
 \right)^2+\left(\frac{\partial F_2}{\partial\zeta}(0,0,\pm 5)
 \right)^2} > 0
$$
and we conclude that, in this example, disease-induced mortality yields 
instability.

The calculations to obtain the conclusion above are standard, but for 
the reader's convenience we give below some details.
With the choices made in this special case we have the following 
kernels where we have highlighted the dependence on $\alpha$
\begin{equation*}
\begin{aligned}
 \Psi_1(\alpha,a)= & \left( 6\Phi(Q^*(\alpha)) - 
         \frac{B^*(\alpha)}{10\Phi(Q^*(\alpha))} \right) \frac{3}{2}
         \sin(2a)\cos(a) \cdot \\
         & \cdot \left( e^{-\alpha a}
        + \alpha\int_0^a e^{-\alpha(a-\sigma)
         -\frac{W^*(\alpha)}{2}\sigma^2} d\sigma \right)\, , \\
 \Psi_2(\alpha,a)= & -\alpha B^*(\alpha) \left( 6\Phi
         (Q^*(\alpha)) - \frac{B^*(\alpha)}{10\Phi(Q^*(\alpha))} \right) 
         \frac{3}{2} \cdot\\
         &\cdot \int_a^\frac{\pi}{2} \sin(2\sigma) \cos(\sigma)
         (\sigma-a) \int_{\sigma-a}^\sigma e^{-\alpha(\sigma-\rho)
         -\frac{W^*(\alpha)}{2}\rho^2}d\rho\, d\sigma\, , \\
 \Psi_3(\alpha,a)= &W^*(\alpha) \frac{3}{2}\sin(2a)\cos(a)\int_0^a
         \sigma e^{-\alpha(a-\sigma)-\frac{W^*(\alpha)}{2}\sigma^2}
         d\sigma \\
 \Psi_4(\alpha,a)= &B^*(\alpha) \frac{3}{2} \int_a^\frac{\pi}{2} 
         \sin(2\sigma)\cos(\sigma)(\sigma-a)\cdot\, , \\
         & \cdot \left( e^{-\frac{W^*(\alpha)}{2}\sigma^2} - \alpha \int_{\sigma-a}^\sigma e^{-\alpha(\sigma-\rho)
           -\frac{W^*(\alpha)}{2}\rho^2} d\rho \right) d\sigma\, .
\end{aligned}
\end{equation*}
Then we evaluate the following derivatives (they are the only derivatives we need)
$$
\begin{aligned}
\frac{\partial\Psi_1}{\partial\alpha}(0,a) &= 
   \sin(2a)\cos(a) \left( - 3W^*(0) + 6a - 6 \int_0^a e^{-\frac{W^*(0)}{2}
   \sigma^2} d\sigma \right), \\
\frac{\partial\Psi_2}{\partial\alpha}(0,a) &= 
   50 \int_{a}^\frac{\pi}{2} \sin(2\sigma)\cos (\sigma) (\sigma-a) 
   \int_{\sigma-a}^\sigma e^{-\frac{W^*(0)}{2}\rho^2} d\rho\, d\sigma\, ,
\end{aligned}
$$
where we have used 
\begin{equation*}
{W^*}'(0) = \left( \frac{12}{125}W^*(0) -1 \right) \frac{\frac{\partial H}
 {\partial\alpha}(0,W^*(0)}{\frac{\partial H}{\partial W}(0,W^*(0))}\, ,
\end{equation*}
obtained by Dini's Theorem applied to (\ref{dini}).
The above expressions are used to compute the following quantities
\begin{equation*}
\begin{aligned}
A & = 1-\int_0^\frac{\pi}{2}\cos(5a)\Psi_4(0,a)\, da \simeq 0.83432,\\
B & = \int_0^\frac{\pi}{2}\sin(5a)\Psi_4(0,a)\, da \simeq 0.236031,\\
C & = \int_0^\frac{\pi}{2}\cos(5a)\frac{\partial\Psi_1}{\partial\alpha}(0,a) da \simeq 2.52632,\\
D & = \int_0^\frac{\pi}{2}\sin(5a)\frac{\partial\Psi_1}{\partial\alpha}(0,a)
da \simeq -0.397545,\\
E & = \int_0^\frac{\pi}{2}\cos(5a)\Psi_3(0,a)\, da \simeq -0.147547,\\
F & = \int_0^\frac{\pi}{2}\sin(5a)\Psi_3(0,a)\, da \simeq -0.177344,\\
\end{aligned}
\end{equation*}
\begin{equation*}
\begin{aligned}
G & = \int_0^\frac{\pi}{2}\cos(5a)\frac{\partial\Psi_2}{\partial\alpha}(0,a) da \simeq -0.0744986,\\
H & = \int_0^\frac{\pi}{2}\sin(5a)\frac{\partial\Psi_2}{\partial\alpha}(0,a) da \simeq -0.249824,\\
\end{aligned}
\end{equation*}
that finally allow us to evaluate
$$
\begin{aligned}
\frac{\partial F_1}{\partial\alpha}(0,0,\pm 5) &= -AC-BD-EG+FH\simeq -1.98062, \\
\frac{\partial F_2}{\partial\alpha}(0,0,\pm 5) &= \pm \left( -BC+AD+FG+EH 
 \right) \simeq \mp 0.877897,\\
\frac{\partial F_1}{\partial\zeta}(0,0,\pm 5) &= \frac{5\pi}{32}A + 
 \frac{2}{3} B \simeq 0.5669, \\
\frac{\partial F_2}{\partial\zeta}(0,0,\pm 5) &= \pm \left( \frac{5\pi}{32} B - 
 \frac{2}{3}A \right) \simeq \mp 0.440352,
\end{aligned}
$$
so that
$\zeta'(0) \simeq 1.42878>0$.

\section{Numerical exploration} \label{sec:stab_an}

In the previous Section we have produced an example showing that the parameter  $\alpha$, representing disease-induced mortality, can actually modify the dynamics of the system and that periodic solutions are possible via Hopf bifurcation.
In order to explore  the model in a systematic way, we now resort to a numerical method that allows to determine the roots of the characteristic equation and follow their displacement as $\alpha$ varies.
The method, proposed in \cite{breda}, provides numerical 
approximations to the rightmost part of the characteristic spectrum 
associated to the model linearized around the equilibrium to be 
investigated. It is indeed well-known that the zero solution of this 
latter is asymptotically stable if and only if all the characteristic 
roots have strictly negative real part.

The numerical scheme developed in \cite{breda} is actually devoted to the stability 
analysis of the scalar Gurtin-MacCamy model \cite{gurtin}, but it can 
be extended straightforwardly to the $m$-dimensional system ($m\geq1$)
\begin{equation}
 \label{gurt}
 \left\{ \begin{aligned}
     &\frac{\partial{\bf P}}{\partial t}+
       \frac{\partial{\bf P}}{\partial a}+
       \mathcal{M}(a,{\bf S}(t)){\bf P}(a,t)=0\, , \\
     &{\bf P}(0,t)=\int_0^{a^\dagger}\mathcal{B}(a,{\bf S}(t))
       {\bf P}(a,t)da\, , \\
     &{\bf S}(t)=\int_0^{a^\dagger}\mathcal{G}(a){\bf P}(a,t)da\, , \\
     &{\bf P}(a,0)={\bf P_{0}}(a), 
 \end{aligned}\right.
\end{equation}
where ${\bf P}:[0,a^\dagger]\times[0,+\infty)\rightarrow\mathbb{R}^{m}$ 
is the $m$-vector of population densities, $\mathcal{M},\mathcal{B}:
[0,a^\dagger]\times\mathbb{R}^{n}\rightarrow\mathbb{R}^{m\times m}$ 
are the matrices of mortality and fertility rates, respectively, and 
${\bf S}:[0,+\infty)\rightarrow\mathbb{R}^{n}$ ($n\geq1$) is the 
$n$-vector of population sizes, i.e. a selection of $n$ homogeneous 
population sub-classes through the weight function $\mathcal{G}:
[0,a^\dagger]\rightarrow\mathbb{R}^{n\times m}$.
\begin{figure}[h]
 \begin{center}
   \includegraphics[width=7.5cm]{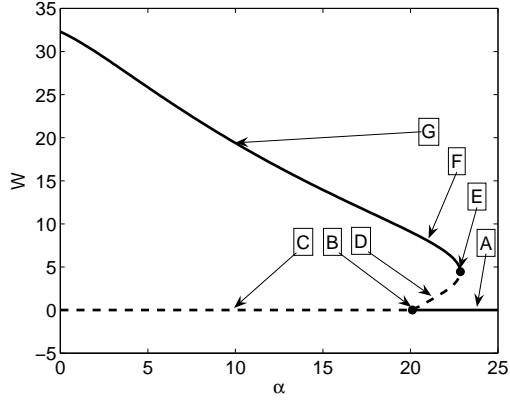}
   \caption{bifurcation diagram of equilibrium $W^{\ast}$ as $\alpha$ 
     varies ($X=34$ in (\ref{plus})).}
   \label{fig3}
 \end{center}
\end{figure}
The epidemic model (\ref{eq}) we are interested in fits into 
(\ref{gurt}) by choosing $m=2$, $n=2$, the population vector
\begin{equation*}
{\bf P}(a,t)=\left(S(a,t),I(a,t)\right)^{T}
\end{equation*}
and the size vector
\begin{equation*}
{\bf S}(t)=\left(Q(t),W(t)\right)^{T} ,
\end{equation*}
which corresponds to selecting
\begin{equation*}
\mathcal{G}(a)=
\begin{pmatrix}
r(a)&r(a)\\
0&q(a)\\
\end{pmatrix}.
\end{equation*}
Moreover, the matrices relative to the vital rates are given by
\begin{equation*}
\mathcal{M}(a,{\bf S}(t))=
\begin{pmatrix}
K(a)W(t)+\mu(a)&0\\
-K(a)W(t)&\alpha+\mu(a)\\
\end{pmatrix}
\end{equation*}
and
\begin{equation*}
\mathcal{B}(a,{\bf S}(t))=
\begin{pmatrix}
R_0^d\beta(a)\Phi(Q(t))&R_0^d\beta(a)\Phi(Q(t))\\
0&0\\
\end{pmatrix}.
\end{equation*}

Within this framework we consider the special choices of Section \ref{secexun},
namely (\ref{choices}) and (\ref{plus}) with $X=34$, letting $\alpha$ vary  
along the corresponding equilibrium curve represented in Figure 
\ref{fig2} (the third curve from the right). 
For each point on the 
curve (included those corresponding to the disease-free equilibrium for 
$W^{\ast}=0$) we are able to compute the rightmost characteristic root 
(rounded to machine precision, see \cite{breda}) and thus to say whether the 
corresponding equilibrium is locally asymptotically stable or not. The 
overall situation is illustrated in Figure \ref{fig3} where 
solid lines denote stability and dashed lines denote instability.

\begin{figure}[!pt]
 \begin{center}
   \includegraphics[width=5.8cm]{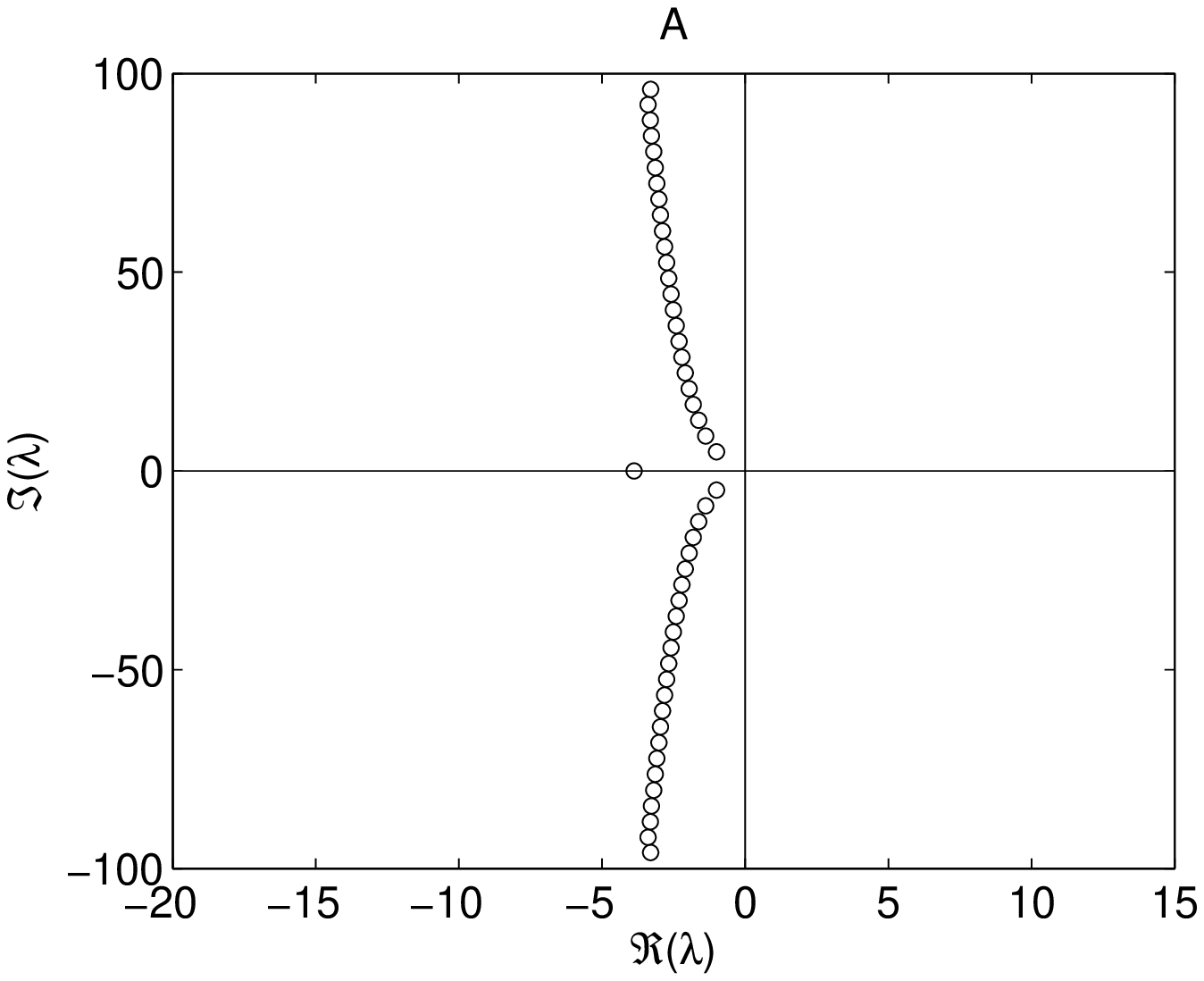}
   \includegraphics[width=5.8cm]{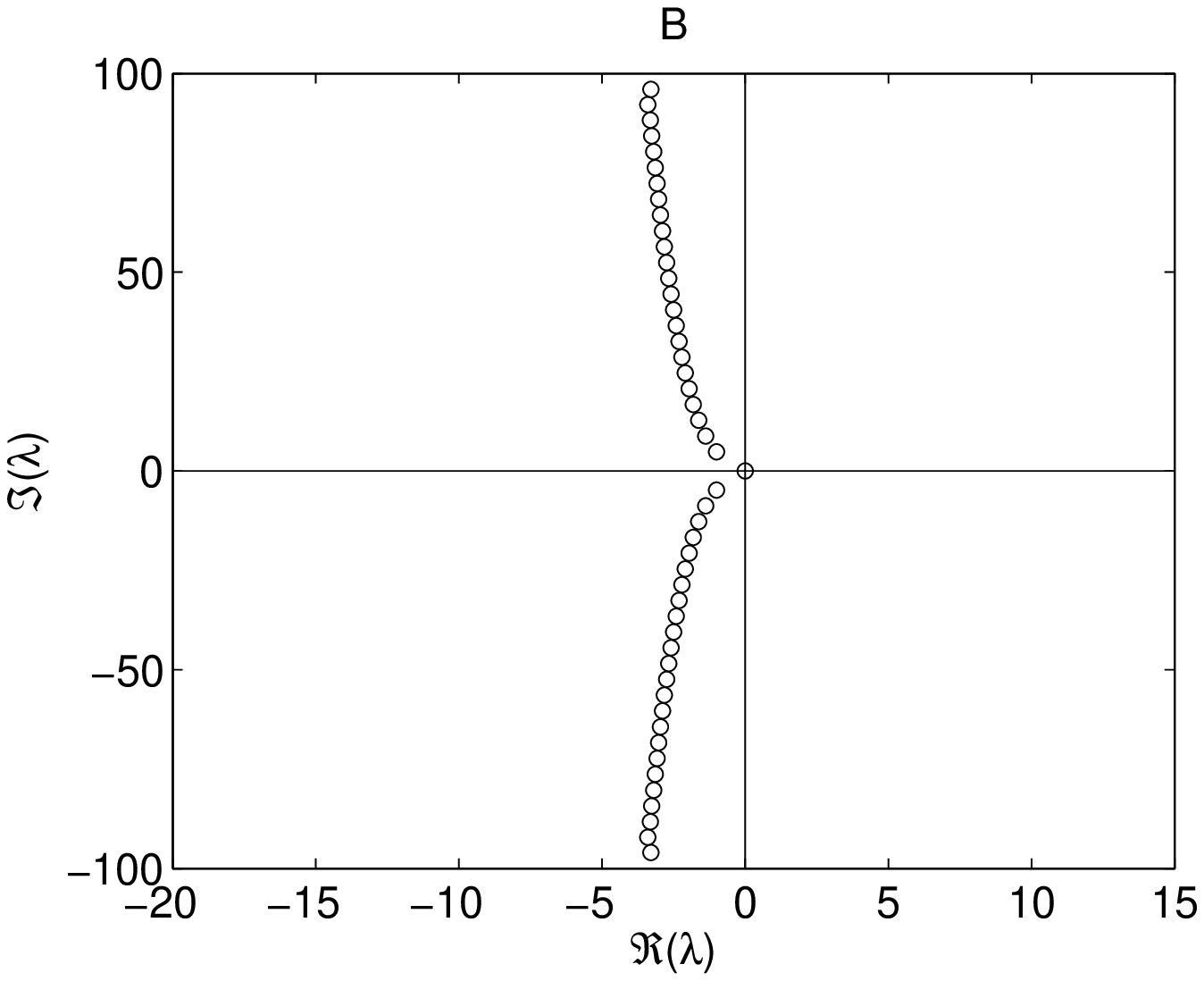}
   \includegraphics[width=5.8cm]{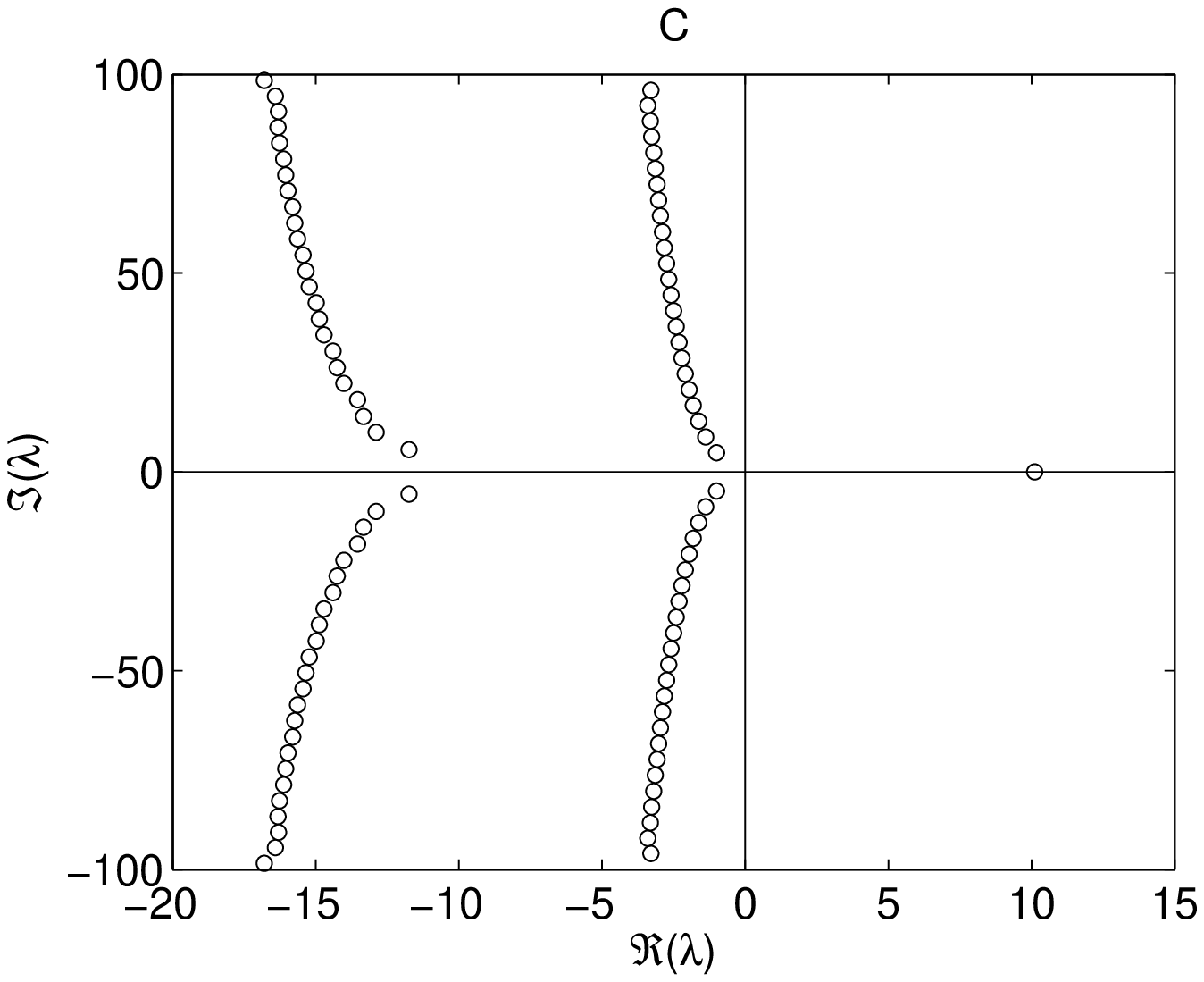}
   \includegraphics[width=5.8cm]{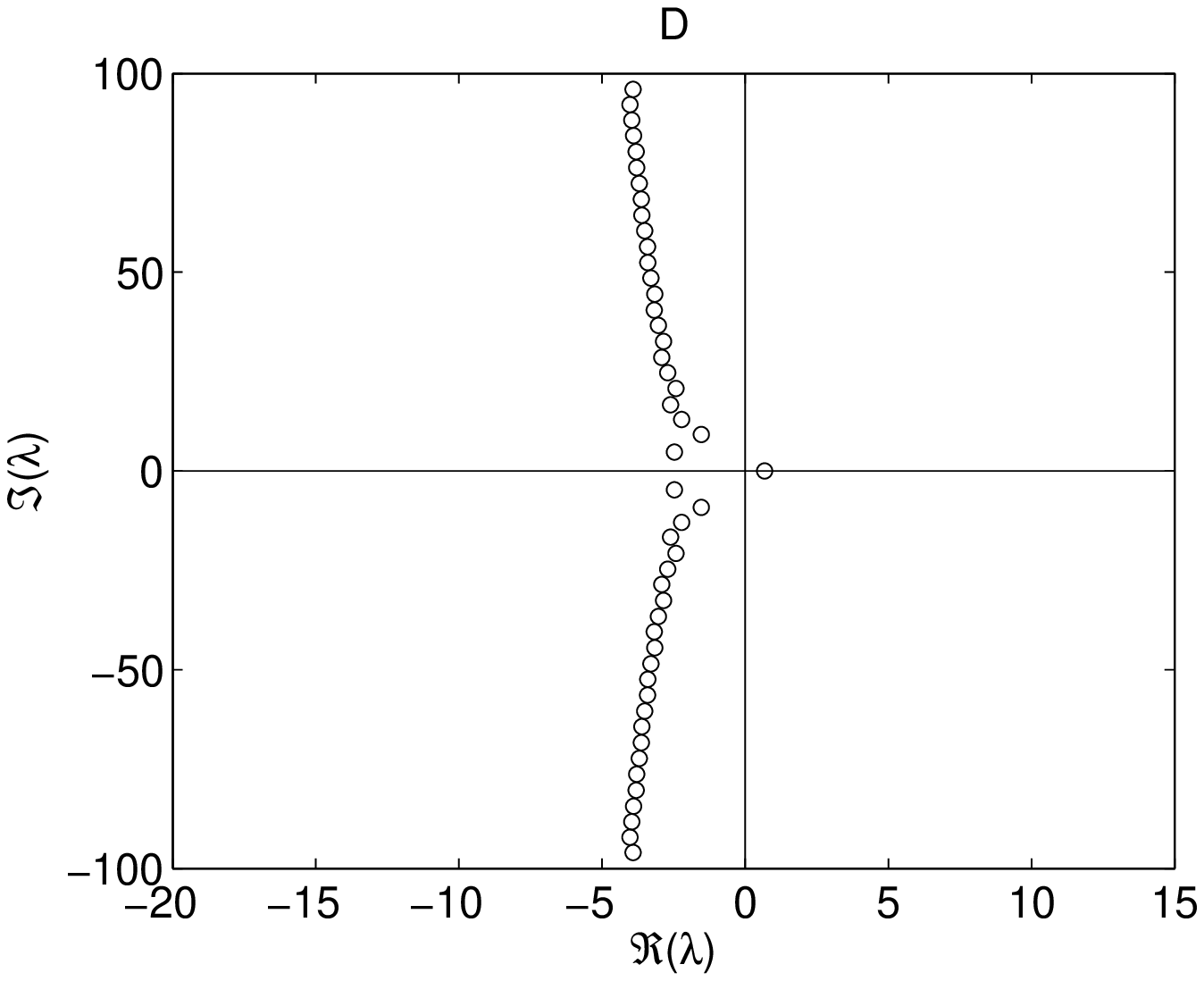}
   \includegraphics[width=5.8cm]{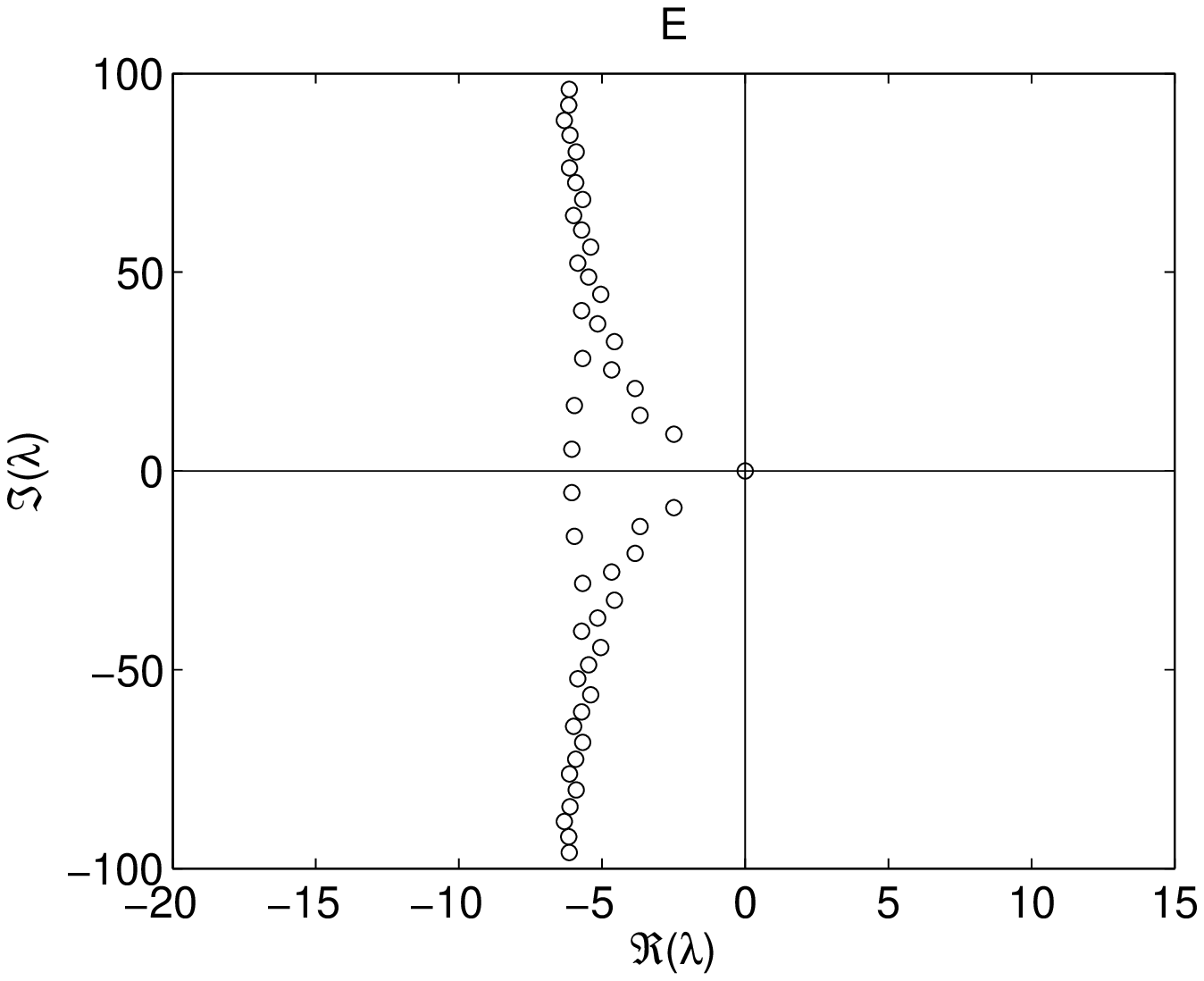}
   \includegraphics[width=5.8cm]{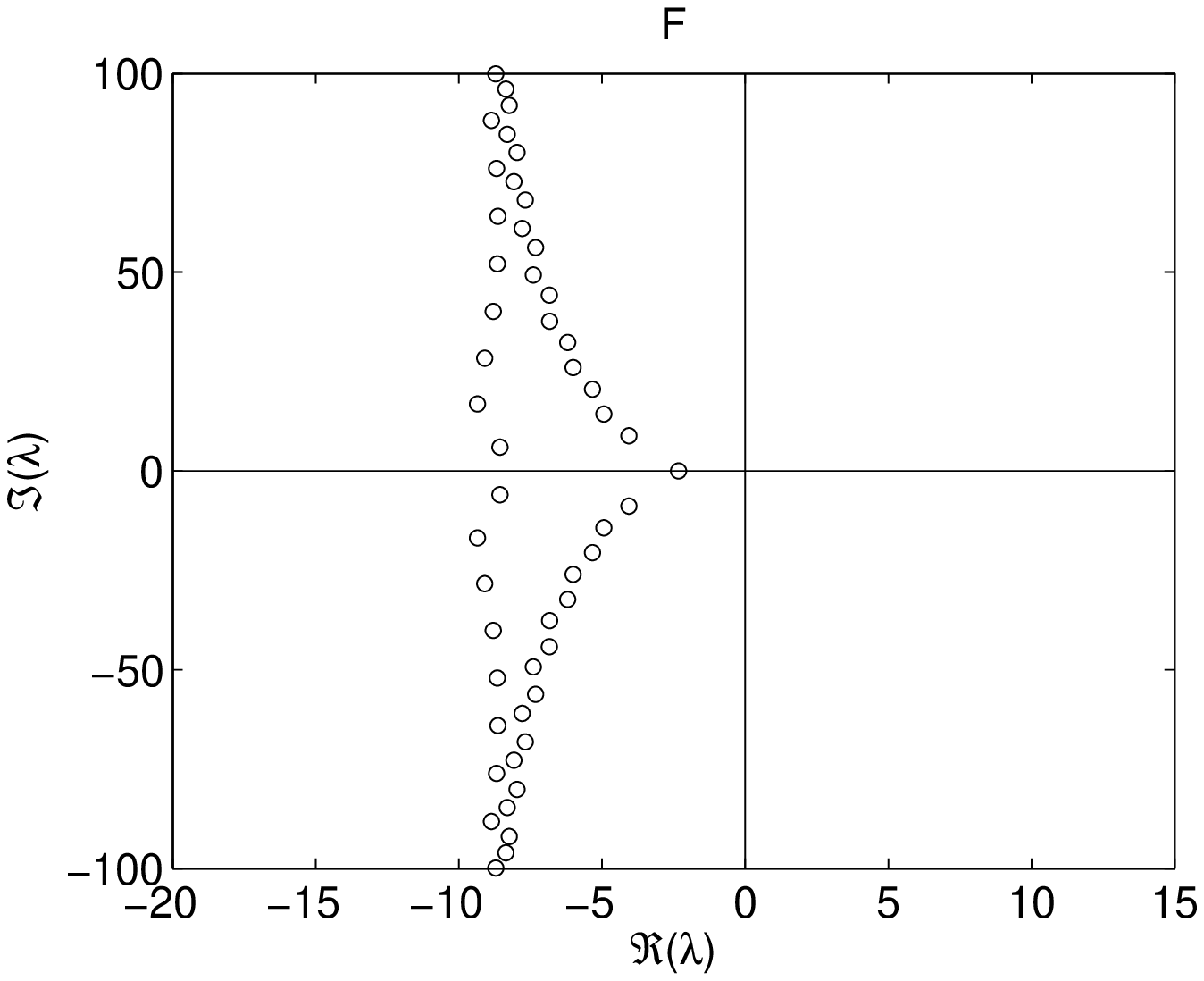}
   \includegraphics[width=5.8cm]{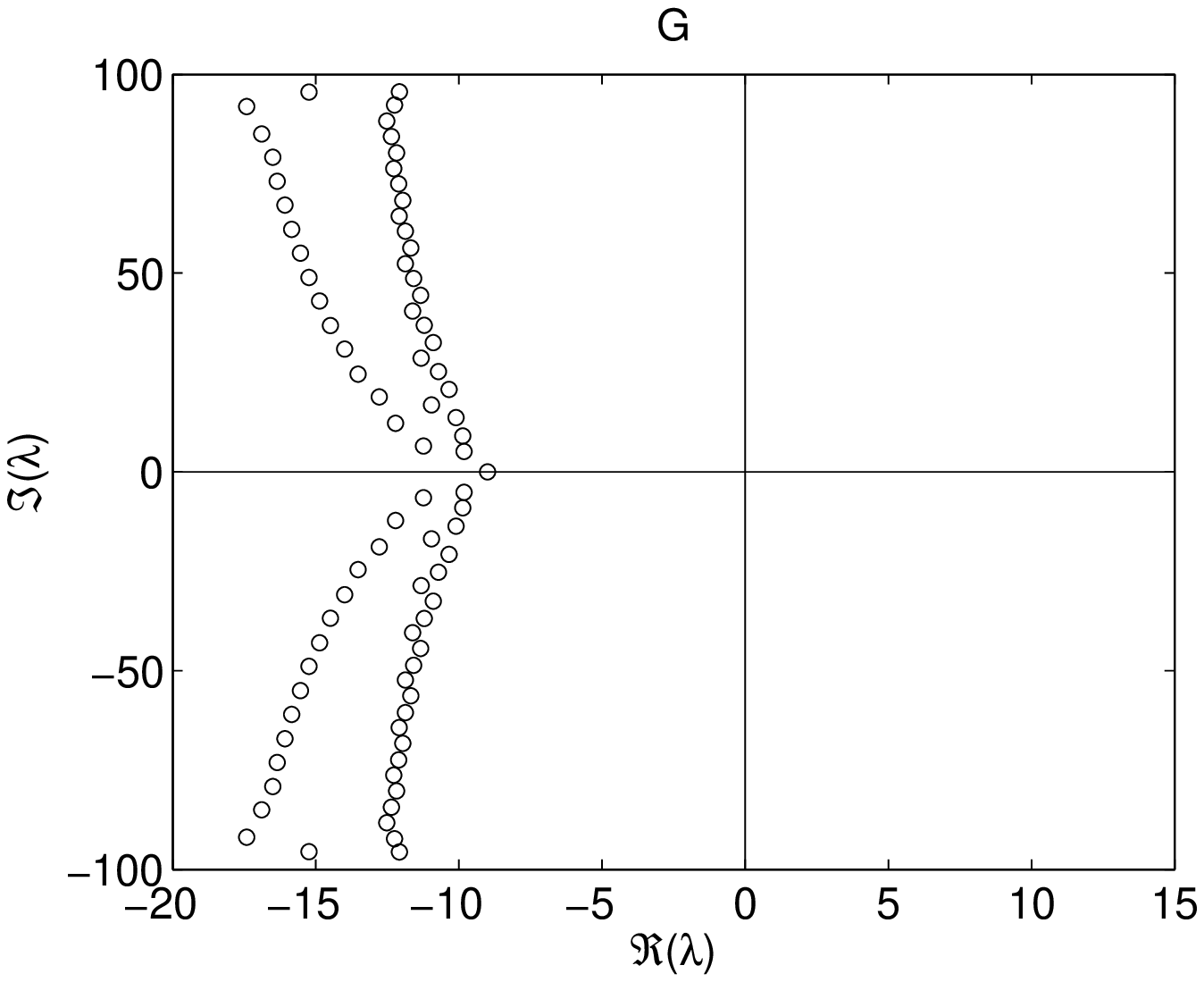}
   \caption{relevant characteristic roots for different 
     equilibrium points, referring to the bifurcation diagram of Figure \ref{fig3} ($X=34$ in (\ref{plus})).}
   \label{fig6}
 \end{center}
\end{figure}

We start our analysis from the right hand branch of the disease free 
equilibrium by investigating, for instance, the spectrum for $\alpha=24$ 
(A in Figure \ref{fig3} and Figure \ref{fig6}): the rightmost roots have 
negative real part 
and hence the trivial equilibrium is stable. By decreasing $\alpha$ the 
real root moves to the right until at $\alpha=\alpha_{1}\simeq20.1143$ 
it crosses the imaginary axis rightward (B). This first bifurcation 
makes the trivial equilibrium lose its stability and a branch of 
endemic equilibria raises. Following the disease-free branch by 
further decreasing $\alpha$ it is confirmed that the instability 
persists (C, $\alpha=10$).

Going back to the bifurcation point at $\alpha_{1}$, we now follow the 
endemic branch by increasing $\alpha$. The branch is unstable (D, 
$\alpha=21$) until at $\alpha=\alpha_{2}\simeq22.8495$ the leading 
root crosses the imaginary axis leftward (E). 
This second bifurcation 
makes the endemic equilibrium gain back its stability. Then the branch 
continues by decreasing $\alpha$ again and at $\alpha=21$ (F) the second 
endemic equilibrium is stable opposite to the first one for the same 
value of $\alpha$ (D). The branch remains stable by further decreasing 
$\alpha$ (F, $\alpha=10$).
\begin{figure}[!htb]
 \begin{center}
   \includegraphics[width=6cm]{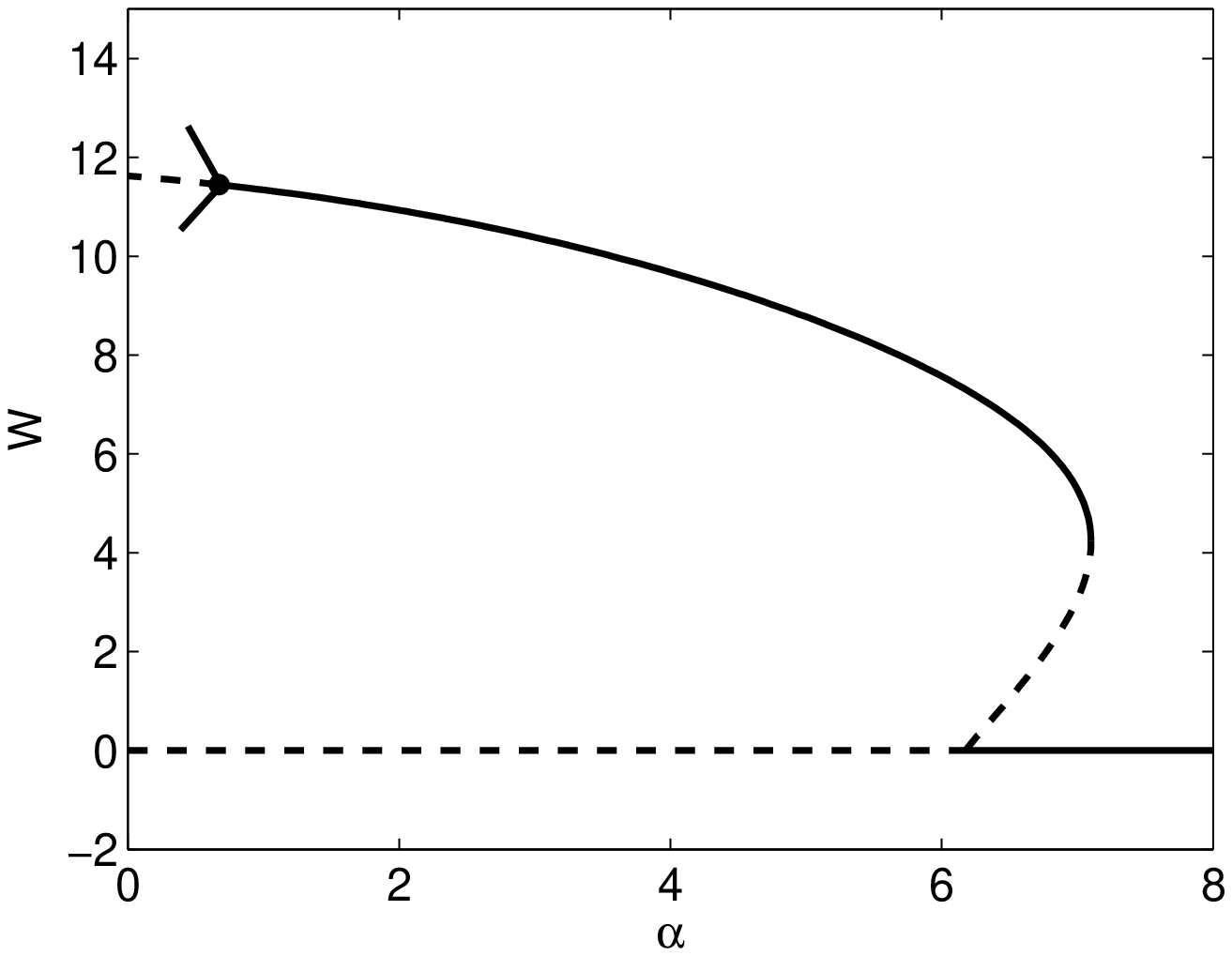}
   \includegraphics[width=6cm]{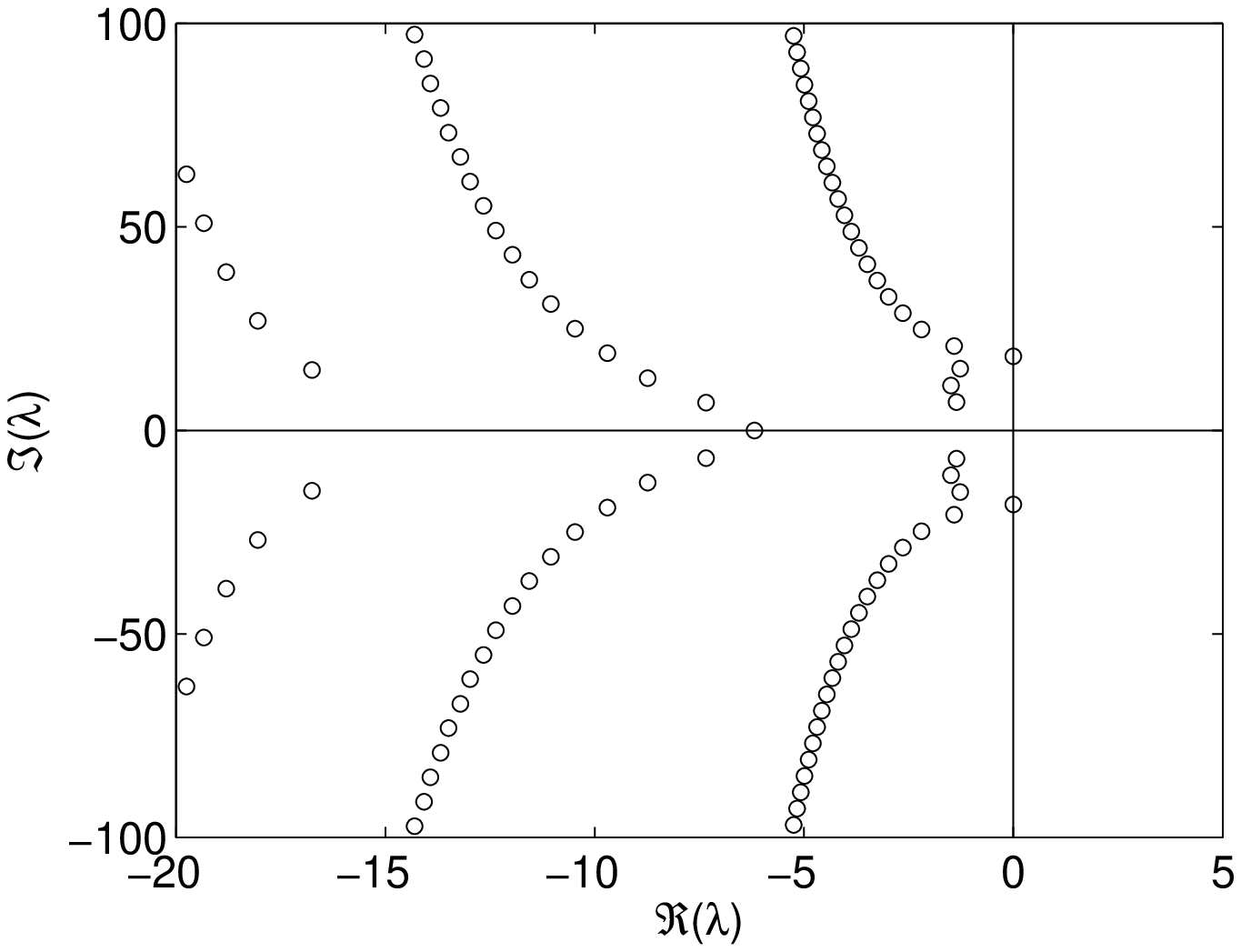}
   \caption{bifurcation diagram of equilibrium $W^{\ast}$ as $\alpha$ 
     varies (left) and relevant characteristic roots (right) for 
     $\alpha\simeq0.6743$, value at which a Hopf bifurcation occurs 
     (black dot in the left figure).}
   \label{fig4}
 \end{center}
\end{figure}

Similar trends are obtained for other values of $X$ for which double 
endemic equilibria exist. On the other side, for those values of $X$ 
for which only one endemic equilibrium exists (for instance $X=10,14$ 
in Figure \ref{fig2}), it can be observed by the roots 
computation that at the first bifurcation the disease-free equilibrium 
loses its stability in favour of the endemic one. This latter 
then preserves its stability for decreasing values of $\alpha$ down to 
$\alpha=0$.

As a consequence of the destabilisation shown in 
Section \ref{secsette}, it is possible that under different choices of the parameters a 
Hopf bifurcation occurs. In fact, 
if we consider the same choices made in (\ref{choices}) for $a^{\dagger}$ and $K$, but set
\begin{equation}
 \label{choices2}
 \begin{array}{l}
   \beta(a)=\frac{3}{2}\sin{(2a)}\, , \\
   q(a)=\sin{(2a)}\, ,\\
   r(a)=\frac{3}{2}\sin{(2a)}\, ,
 \end{array}
 \hspace{2cm}
 \begin{array}{l}
   \mu(a)=\tan a\, , \\
       \Phi(x)=\max\{1-\frac{x}{18},0\}\, ,
 \end{array}
\end{equation}
we get the bifurcation diagram represented in 
Figure \ref{fig4} (left). The black dot, corresponding to $\alpha
\simeq 0.6743$, indicates a Hopf bifurcation through which the 
equilibrium on the stable branch (as usual in solid line) loses 
its stability and a limit cycle arises. The right figure shows the 
existence of the corresponding couple of characteristic roots crossing 
the imaginary axis from left to right as $\alpha$ decreases. Thus, in 
this case the disease-induced mortality has a stabilising effect.

\begin{figure}[!htb]
 \begin{center}
   \includegraphics[width=8cm]{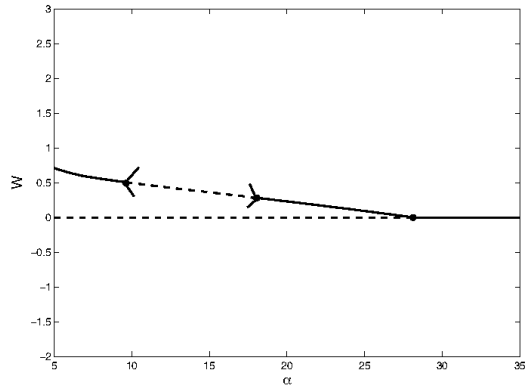}
   \caption{bifurcation diagram of equilibrium $W^{\ast}$ as $\alpha$ 
     varies.  The dotted line corresponds to periodic solutions.}
   \label{fig5}
 \end{center}
\end{figure}

We can also show a choice of the parameters for which the 
diagram is closed and there is, for increasing $\alpha$, a 
destabilisation and, successively, a stabilisation.   In fact, for
\begin{equation}
\label{choices3}
q(a)=10a\, ,\quad r(a)=\frac{3}{5}\sin(2a)\, ,\quad X=15 \quad
 \mbox{and} \quad R_0^d=1.35\, ,
\end{equation}
one has the bifurcation diagram of Figure \ref{fig5}.
%
\section{Conclusions}

The epidemic model studied here differs from the classical ones for 
the population size, which is not constant.   As pointed out in 
\cite{huang}, this leads to the possibility of having multiple endemic 
equilibria.  In our case, which is regulated by a very different 
mechanism, we are able to show that under some conditions 
two endemic equilibria occur. One of them is unstable and the other 
is stable.

It is well known that, when the population size is fixed and 
transmission is inter-cohort, one has 
uniqueness of endemic states.  See for example \cite{bci} and 
\cite{greenhalgh}.  The stability of the steady states is studied in 
\cite{thieme}, \cite{CIM2} and we obtain here an 
analogous stability change for the trivial equilibrium.

For the sake of simplicity we have limited ourselves to investigating 
the case of inter-cohort transmission.  Allowing for a general 
transmission kernel $\lambda(a,t)=\int_0^{a^\dagger}K(a,\sigma)
I(\sigma,t)\, d\sigma$ makes the analysis of the model very difficult 
already in the case of constant population size (see \cite{Inaba}), 
while we wished to emphasise the peculiarities due to 
infection-related deaths.

The situation where multiple endemic equilibria exist may not 
be realistic.  In fact, the model studied assumes that contagion 
can happen only in the first and in the last period of an individual 
life.   A somehow similar situation was studied in \cite{busenberg}, 
where the infection rate is piecewise constant.   More precisely, 
juveniles and noncore adults cannot be infected, while only core 
adults can.

We have shown that the model with age-structured contact rates and 
variable population, because of infection-related deaths, has a much 
richer bifurcation diagram than models with only one of these features.
For example, if (\ref{choices-stab}), (\ref{hyp-radim}) 
and (\ref{choices-stab2}) hold, there is a change of stability and the 
presence of the extra-mortality destabilises the endemic equilibrium.   
Moreover, if either (\ref{choices2}) or (\ref{choices3}) holds, Hopf 
bifurcations occur.

We believe to have shown the main 
theoretical steps necessary to determine the stationary solutions of 
the model and their stability properties, adding the use of the 
numerical tool provided in \cite{breda} whenever analytical conclusions 
are hard, when not impossible, to draw.   We also have highlighted the 
complexity of this model, that presents different kinds of dynamics.

\section*{Acknowledgements}

The authors wish to thank M. Iannelli and A. Pugliese for having shared 
their priceless knowledge of the field and for their fruitful 
observations on the proposed endemic model.


\end{document}